\documentclass[3p,times]{elsarticle}

\usepackage{color}
\usepackage{slashbox}
\usepackage{dsfont}
\usepackage{amsmath}
\usepackage{amssymb}
\usepackage{graphicx}
\usepackage{enumitem} 
\usepackage{breqn}
\usepackage{amsthm}
\theoremstyle{definition}
\newtheorem{Theorem}{Theorem}[section]
\newtheorem{Proposition}{Proposition}[section]
\newtheorem{Lemma}{Lemma}[section]

\newtheorem{Example}{Example}

\newenvironment{Proof}{\paragraph{Proof:}}{\hfill$\blacksquare$}
\newcommand{\norm}[1]{\left\lVert#1\right\rVert}
\newcommand{\E}{\mathbb{E}}

\usepackage{amssymb}
\usepackage{subfigure}

\usepackage[figuresright]{rotating}

\begin{document}
	
	\begin{frontmatter}
		
		
		
		

		\title{A tempered subdiffusive Black-Scholes model}
		\author[label1]{Grzegorz Krzy\.zanowski}
		\ead{grzegorz.krzyzanowski@pwr.edu.pl}
		
		\author[label1]{Marcin Magdziarz}
		\ead{marcin.magdziarz@pwr.wroc.pl}
		
		
		\address[label1]{Hugo Steinhaus Center,
			Faculty of Pure and Applied Mathematics, Wroclaw University of Science and Technology
			50-370 Wroclaw, Poland}
				\begin{abstract} 
			In this paper, we focus on the tempered subdiffusive Black-Scholes model. The main part of our work consists of the finite difference method as a numerical approach to the option pricing in the considered model. We derive the governing fractional differential
			equation and the related weighted numerical scheme.
			The proposed method has the $2-\alpha$ order of accuracy with respect to time, where $\alpha\in(0,1)$ is the subdiffusion parameter, and $2$ with respect to space.
			Furthermore, we provide the stability and convergence analysis.
			Finally, we present some numerical results. 
			
		\end{abstract}%
		%
		
	\end{frontmatter}
	\section*{Introduction}
Options are one of the most popular and important financial derivatives, therefore the question about their valuation
has an essential meaning for financial institutions and global economies. The value of the global derivatives market is estimated to $700$ trillion dollars to upwards of $1,5$ quadrillion dollars \cite{market}. In $2019$ the volume of traded derivative contracts reached $34,47$ billion, including $15,23$ billion of options contracts \cite{fia}. 
The option is a contract in which the holder can buy/sell a property for a fixed price $K$. 
Underlying assets can be actions, stock exchange indexes, foreign currency, futures contracts, or obligations.

Over the past two decades, the B-S model has been increasingly attracting interest as an effective and easy tool for option valuation. The model was of such great importance that
the authors were awarded the Nobel Prize for Economics in 1997. Although the discovery was initially recognized as outbreaking, the model
can not be used in many different cases \cite{ladde2009development}.
It was the reason the model was generalized for allowing such features as, e.g., stochastic interest or volatility \cite{angelini2014delta, goncu2014efficient, malhotra2022pricing}, transaction costs \cite{arregui2012numerical, de2021pseudospectral, yan2021utility}, jumps \cite{boen2020european, costabile2014option, wang2021efficient}, and switching regime \cite{costabile2014option, elliott2014pricing, koleva2017numerical}.

\begin{figure}[h]
	\centering
	\includegraphics[scale=0.6]{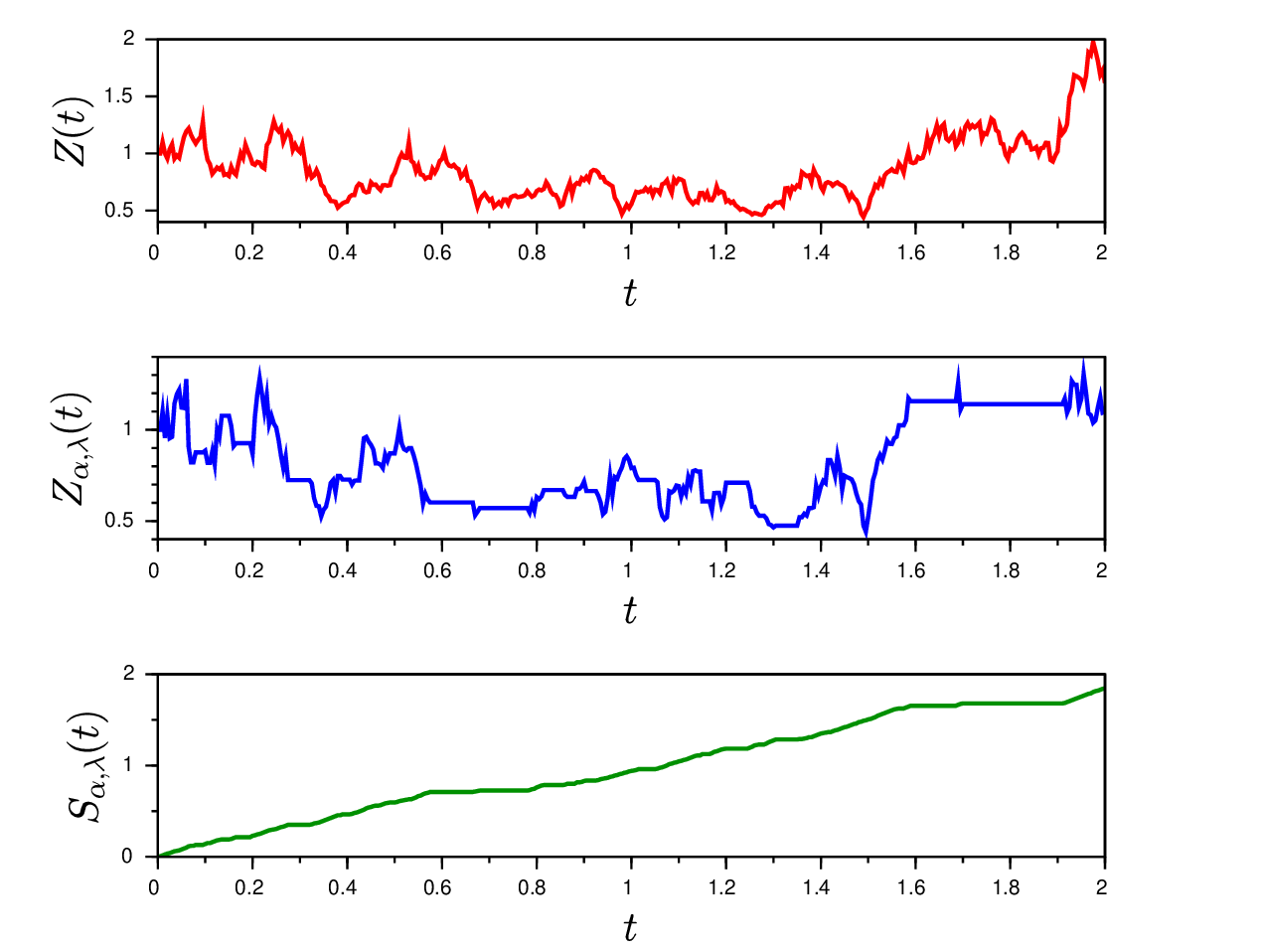}
	\caption{The sample trajectory of GBM (up) with its tempered subdiffusive analogue (middle) and the corresponding inverse subordinator (down). In the tempered subdiffusive GBM, the constant periods characteristic for
		emerging markets can be observed. The parameters are $Z_{0}=\sigma=\mu=\lambda=1$, $\alpha=0.9$.}\label{Figure1}  
\end{figure}

In recent years, it can be observed among different economies characteristic periods of stagnation (see, e.g. \cite{bor, janczura},
and the references in them). This feature is most common for
emerging markets in which the number of participants and thus the number of transactions,
is rather low. These characteristic periods of financial processes correspond to the trapping events in which
the test particle is motionless \cite{eliazar}. In response to empirical evidence of fat tails, $\alpha$-stable distribution as an alternative to the Gaussian 
law was proposed. The subdiffusive regime is obtained by the use of the inverse stable
subordinator (see \cite{Gajda, orzel} and references in them).
The stable distribution has found many important applications, for
example in physics \cite{16,15,17} and electrical engineering \cite{18}. 
Since Mandelbrot \cite{mandelbrot} and Fama \cite{fama} introduced the $\alpha$-stable distribution in modeling financial
asset returns, numerous empirical studies have been conducted in both natural and
economic sciences. In \cite{rachev1, rachev2} and the references therein, a wide range of applications of the $\alpha$-stable distribution in finance
is considered. The empirical study confirms that the $\alpha$-stable distribution does not always adequately describe the data following the stagnation phenomenon. 
Asset return time series often demonstrate heavier tails than the normal distribution and thinner tails than the $\alpha$-stable distribution.
Partly in response to the above empirical evidences, and to
maintain suitable properties of the stable model, a proper generalization of the
$\alpha$-stable distribution was introduced. The $\alpha$-stable $\lambda$-tempered distributions are manifested by heavier tails than the
normal distribution and thinner than the stable distribution, moreover, they have finite moments of all orders \cite{rosinski}.
Similarly to subdiffusive B-S, tempered subdiffusive Black-Scholes (tsB-S) is the generalization of the classical B-S model to the cases where the underlying assets
display characteristic periods in which they remain motionless \cite{ja3}. The standard B-S model assumes that the asset is described by a continuous Gaussian random walk, so the underlying asset at each step has to move up or down. As a result of option pricing in such a stagnated market, the fair price provided by the B-S model is misestimated.
To properly describe this dynamics, the tsB-S model assumes that the underlying asset is driven by $\alpha$-stable $\lambda$-tempered inverse subordinator (see \cite{Gajda, orzel} and the references therein). The frequency of the constant periods appearing then depends on the subdiffusion parameter $\alpha\in(0,1)$ and tempering parameter $\lambda>0$ (the particular case of tsB-S where $\lambda=0$ is the subdiffusive B-S considered e.g. in \cite{ja2, MM, zhang}). If $\alpha\rightarrow1$, tsB-S is reduced to the classical model.
In contrast to the tsB-S model, the B-S model does not take into account the empirical property of the stagnated periods of the underlying instrument. 
In Figure \ref{Figure1}   we compare the sample simulation of the underlying instrument in the classical and tempered stable market model. Even a short constant period of a market cannot be simulated by the classical B-S model.
As a generalization of the classical B-S model, the tsB-S model can be used in a wide range of cases, including all cases where B-S can be applied. 

In this paper, we find the corresponding
fractional differential equation, and for such a model we solve the problem numerically. To do so, we use the weighted finite difference method. We provide the stability/convergence analysis. Finally, we present some numerical examples.
The most important advantages of our model are the property of the lack of arbitrage and the clear motivation of such a generalization of the
classical and subdiffusive B-S model. 
\section{Tempered subdiffusive B-S model} 
\subsection{Assumptions of the tsB-S}
Let us consider a market whose evolution is occurring up to the time horizon $T$ and is contained in the probability space 
$(\Omega, \mathcal{F} , \mathbb{P})$.
Here, $\Omega$ is the sample space,  $\mathcal{F}$ is filtration interpreted as information on the history of the asset price which completely is available for the
investor and $\mathbb{P}$ is the "objective" probability measure. The assumptions are the same as in the classical case \cite{ja} with the exception that we do not have to assume market liquidity and that 
the underlying instrument instead of the Geometric Brownian Motion (GBM) has to follow a tempered subdiffusive GBM \cite{Gajda}:
$$
\left\{ \begin{array}{ll}
	Z_{\alpha,\lambda}\left(t\right)=Z\left(S_{\alpha,\lambda}(t)\right),\\
	Z\left(0\right)=Z_{0},                                                         
\end{array}\right.$$
where $0<\alpha<1$, $\lambda>0$, $Z_{\alpha,\lambda}\left(t\right)$ is the price of the underlying instrument, $\mu$ - drift (constant), $\sigma$ - volatility (constant), $B(t)$- Brownian motion, $Z(t)=Z\left(0\right)e^{\mu t+\sigma B(t)}$, $S_{\alpha,\lambda}\left(t\right)$ is the inverse $\alpha$-stable $\lambda$ -tempered subordinator defined as 
$S_{\alpha,\lambda}\left(t\right)=\inf\left(\tau:W_{\alpha,\lambda}(\tau)>t\right)$ \cite{ken1999levy}, $W_{\alpha,\lambda}(\tau)$ is a $\alpha$-stable $\lambda$ -tempered subordinator defined by its Laplace transform
$\E e^{-uW_{\alpha,\lambda}(t)}=e^{-t((u+\lambda)^{\alpha}-\lambda^{\alpha})}$ \cite{Meser}. We assume that $S_{\alpha,\lambda}(t)$ is independent of $B(t)$ for each $t\in[0,T]$.

Note that with $\alpha\rightarrow 1$ the tsB-S model reduces to the classical case. Due to its simplicity and practicality, the classical B-S model
is one of the most widely used in option pricing. Although in contrast to the subdiffusive and tempered subdiffusive cases, it does not take into account the empirical property of constant price periods.
The method of calibrating $\alpha$ and $\lambda$ from empirical data is the same as in \cite{Wyl}.

\newpage

As shown in \cite{jaCRR} for tsB-S the put-call parity holds:
\begin{Proposition}
	For the fair price of the European call and put options in tsB-S, we have the following relationship:
	\begin{equation}
		C_{tsB-S}=P_{tsB-S}+Z_{0} -\mathbb{E} Ke^{-rS_{\alpha,\lambda}(T)},\label{putcallpar}
	\end{equation}
\end{Proposition}
where $C_{tsB-S}$ and $P_{tsB-S}$ are the fair prices of the European put and call options in the tsB-S model. Note that both $C_{tsB-S}$ and $P_{tsB-S}$ should have identical $Z_{0}$, $K$, $T$,$\sigma$, $r$, $\alpha$, $\lambda$.
Here and in the entire paper: $K$ - strike,
$T$ - expiration time,
$\sigma$ - volatility,
$r$ - interest rate. Without loss of generality, in the entire paper we assume that the dividend rate is equal $0$. One of the most expected properties of the market is that there is no possibility of earning money without taking the risk. This property is called 
the lack of arbitrage and formally means that the self-financing strategy $\phi$ that leads to a positive profit without any probability of intermediate loss can not be constructed \cite{4}.
According to the Fundamental Theorem of Asset Pricing \cite{4}, the market model described by $(\Omega, \mathcal{F}, \mathbb{P})$ and
underlying instrument $Z_{\alpha,\lambda}$ with filtration $\mathcal{F}_{t\in[0,T]}$ is arbitrage-free if and only if there exists
a probability measure $\mathbb{Q}$, (called the risk neutral measure) equivalent to $\mathbb{P}$ such that the asset $Z_{\alpha,\lambda}$ is a martingale
with respect to $\mathbb{Q}$. Under this measure,
financial instruments have the same expected rate of return, regardless of the variability of the prices.
This contrasts with the physical probability measure (the actual probability
distribution of prices), under which more risky instruments have a higher expected rate of return
than less risky instruments. Let us introduce the probability measure \begin{equation}
	\mathbb{Q}(A)=\int_{A}e^{-\gamma B(S_{\alpha,\lambda}(T))-\frac{\gamma^{2}}{2}S_{\alpha,\lambda}(T)}dP,\label{1}
\end{equation}
where $\gamma=\frac{\mu+\frac{\sigma^2}{2}}{\sigma}$, $A\in\mathcal{F}$.
As shown in \cite{Gajda} the process $Z_{\alpha,\lambda}$ is a martingale with respect to $\mathbb{Q}$, so we have the following \newpage

\begin{Theorem}\cite{Gajda}
	The tsB-S model is arbitrage-free.
\end{Theorem}

Another property of the market model is the so-called completeness. Intuitively, the market model is complete
if the set of possible gambles on future states of the world can
be constructed with existing assets. More formally, the market model is complete if every
$\mathcal{F}_{t\in[0,T]}$ -measurable random variable $X$ admits a replicating self-
financing strategy $\phi$ \cite{4}.

The Second Fundamental Theorem of Asset Pricing \cite{4} states that a
market model described by $(\Omega, \mathcal{F}, \mathbb{P})$ and
underlying instrument $Z_{\alpha,\lambda}$ with filtration $\mathcal{F}_{t\in[0,T]}$ is complete if and only if there is a unique
martingale measure equivalent to $\mathbb{P}$. 
\begin{Theorem}\cite{Gajda}
	The market model in which the price of the underlying instrument follows the tempered subdiffusive GBM $Z_{\alpha,\lambda}$
	is incomplete.
\end{Theorem}

Market incompleteness means that there is no unique fair
price of financial derivatives because for different martingale measures, different prices could be obtained. 
Although $\mathbb{Q}$ defined in (\ref{1}) is not unique, in the sense of the criterion of minimal relative entropy, it is the ``best'' martingale
measure. 
It means that the measure $\mathbb{Q}$ minimizes the distance to the
measure $\mathbb{P}$ \cite{Gajda}. Another essential fact is that for $\alpha\to1$, $\mathbb{Q}$ reduces to the measure of the classical
B-S model which is arbitrage-free and complete. It is consistent with our intuition if we consider the tsB-S model as a generalization of the standard B-S model.
Therefore, in this paper we will use the martingale measure $\mathbb{Q}$ defined in (\ref{1}) as a reference measure.

\subsection{The fair price of a call option in the tsB-S model}

In this section, we will prove a tempered analogue of Theorem 2.3. from \cite{ja2}.
\begin{Theorem}\label{TW1q}
	Let us denote $f(x)\sim g(x)$ with $x\to\infty$ if $\lim_{x\to\infty}\frac{f(x)}{g(x)}=1$. We introduce the following variable:
	\begin{equation}
		x=\ln z \label{10}     
	\end{equation}	
	and function:
	\begin{equation}
		u\left(x,t \right)=v\left(e^{x},T-t \right).\label{11}     
	\end{equation}	
	Then the fair price of a call option in the tsB-S model with respect to $\mathbb{Q}$ is equal to $v(z,0)$,
	where $v(z,t)$ satisfies (\ref{10}) and (\ref{11}), and $u(x,t)$ is 
	the solution of
	\begin{equation}
		\begin{cases}
			\partial^{\alpha,\lambda}_{t} u \left(z,t\right)= \displaystyle\frac{1}{2}\sigma^2\displaystyle\frac{\partial^2 u\left(x,t\right)}{\partial x^2}+\left(r-\displaystyle\frac{1}{2}\sigma^2\right)\displaystyle\frac{\partial u\left(x,t\right)}{\partial x}-ru\left(x,t\right),\\
			u\left(x,0\right)=\max\left(e^x-K,0\right),\\
			\lim\limits_{x\rightarrow -\infty} u\left(x,t\right)=0,\\
			u\left(x,t\right)\sim e^{x} \quad\text{for}\quad x \rightarrow \infty,\\
		\end{cases}\label{12}     
	\end{equation}
	for $\left(x,t\right)\in\left(-\infty,\infty\right)\times(0,T]$. The operator $\partial^{\alpha,\lambda}_{t}$ is a tempered Caputo fractional derivative defined as \cite{Meser}:
	\begin{equation*}\partial^{\alpha,\lambda}_{t} u(z,t)=e^{-\lambda t} \text{ }_{RL}D_{t}^{\alpha} (e^{\lambda t} (u(z,t)-u(z,0)))-\lambda^{\alpha}(u(z,t)-u(z,0)),
	\end{equation*}
	where $\alpha\in(0,1)$, $\lambda>0$, and the Riemann-Louville derivative $\text{ }_{RL}D_{t}^{\alpha}$ is defined as \cite{kilbas2006theory}:
	\begin{equation*}
		_{RL}D_{t}^{\alpha}g(t)=\frac{1}{\Gamma(1-\alpha)}\frac{d}{dt}\int_{0}^{t}\frac{g(s)ds}{(t-s)^{\alpha}}
	\end{equation*}
	for locally integrable function $g$ on $(0,T)$.
\end{Theorem}
\begin{Proof}
	\label{49o}
		Let us consider a fair price for the European call option in the standard B-S model $\xi(z,t)$ depending on the price of the underlying instrument $z=Z_{0}$ and the time $t$ left to expiration $T$, for $(z,t)\in(0,\infty)\times[0,T]$. In other words, $\xi(z,T-t)$ is the fair price of this option depending on the actual time $T-t$. By the B-S formula \cite{wilmott1995mathematics}, it holds that: \[\xi(z,t)=z\Phi(d_{1})-Ke^{-rt}\Phi(d_{2}),\] where \[d_{1}=\frac{\log(\frac{z}{K})+(r+\frac{1}{2}\sigma^{2})t}{\sigma\sqrt{t}},\]
		\[d_{2}=\frac{\log(\frac{z}{K})+(r-\frac{1}{2}\sigma^{2})t}{\sigma\sqrt{t}},\]
		and $\Phi$ is is the CDF of the normal distribution. For $(z,t)\in(0,\infty)\times(0,\infty)$ $\xi(z,t)$ follows a classical B-S equation \cite{ja, McDonald,wilmott1995mathematics}:
		\[\partial_{t}\xi(z,t)= \displaystyle\frac{1}{2}\sigma^2 z^2 \displaystyle\frac{\partial^2 \xi\left(z,t\right)}{\partial z^2}+rz\displaystyle\frac{\partial \xi\left(z,t\right)}{\partial x}-r\xi\left(z,t\right).\]
		Furthermore, there is an initial condition for a call option $\xi(z,0)=\max(z-K,0)$, for $z\in(0,\infty)$.
		Let us consider $w(z,t)$ given by
		\begin{equation}
			\displaystyle w(z,t)=\int_{0}^{\infty} \xi(z,s) \vartheta_{\alpha,\lambda}(s,t) ds,\label{tylkomk}
		\end{equation}
		where $(z,t)\in(0,\infty)\times(0,\infty)$ and $\vartheta_{\alpha,\lambda}(s,t)$ denotes a density of $S_{\alpha,\lambda}(t)$. Based on \cite{Meser}, formula (6.3), we know that for $(z,t)\in(0,\infty)\times(0,T]$, $w(z,t)$ follows
		\begin{equation}\partial^{\alpha,\lambda}_{t}w(z,t)= \displaystyle\frac{1}{2}\sigma^2 z^2 \displaystyle\frac{\partial^2 w\left(z,t\right)}{\partial z^2}+rz\displaystyle\frac{\partial w\left(z,t\right)}{\partial x}-rw\left(z,t\right)\label{4lataa}\end{equation}
		and \begin{equation}
			w(z,0)=\xi(z,0).\label{nnnnn}
	\end{equation}
We introduce the function $v(z,t)=w(z,T-t)$ for $(z,t)\in(0,\infty)\times[0,T]$. Then, we have: 
		\begin{equation}
			v(z,0)=w(z,T)=\int_{0}^{\infty}\xi(z,s)\vartheta_{\alpha,\lambda}(s,T)ds.\label{ssqq}
		\end{equation}
		Based on \cite{stanislavsky2003black,Gajda, MM} the right side of (\ref{ssqq}) is a fair price (in respect with $\mathbb{Q}$) of European option for the same payoff as $\xi$. Therefore, we conclude that $v(z,0)$ is a fair price of the European call option.
		Let us introduce a variable $x=\ln z$ and a function $\Upsilon(x,t)=\xi(e^{x},t)$, $u(x,t)=w(e^{x},t)$ for $(z,t)\in(0,\infty)\times[0,\infty)$. Based on (\ref{4lataa}) and (\ref{nnnnn}) $u(x,t)$ follows two first conditions of (\ref{12}). By (\ref{tylkomk}) and the definition of the function $\Upsilon$ and $u$ we get 
		\begin{equation}
			\displaystyle u(x,t)=\int_{0}^{\infty} \Upsilon(x,s) \vartheta_{\alpha,\lambda}(s,t) ds,\label{tylko}
		\end{equation}	
		where $(x,t)\in(-\infty,\infty)\times[0,T]$.
		Let us observe that, based on (\ref{tylko}), the boundary conditions of $\Upsilon$ (i.e., the boundary conditions of the standard B-S equation after space transformation) will be conserved for $u$, therefore we obtain the last two conditions of (\ref{12}). Let us observe that the boundary conditions are related with the financial interpetation - i.e., for sufficiently low price of the underlying instrument the option is useless, and that for high enough price of the underlying we use the option obtaining the value of this asset minus the discounted value of $K$.
		It is important to note that \cite{Meser}, by (\ref{tylko}) we obtain the existence and uniqueness of (\ref{12}). Moreover, the price of the European call option in repsect with $t$ in the classical B-S model is bounded and $\vartheta_{\alpha,\lambda}(s,t)$ is a probability density of $S_{\alpha,\lambda}(t)$. Therefore, by (\ref{tylko}) we find that $u(x,t)$ for $(x,t)\in(-\infty,\infty)\times(0,T]$ is finite.
\end{Proof}\\
Please note, that a proof of Theorem \ref{TW1q} for $\lambda=0$ serves as a proof of Theorem 2.3. from \cite{ja2}.
Note that for $\eta(t)=\xi(z,t)-\xi(z,0)$ we have $\eta(0)=0$. We recall that for all differentiable functions in the sense of the Caputo and Riemann-Louville, following the zero initial
condition property, the Riemann-Louville derivative is equal to its Caputo equivalent \cite{kilbas2006theory}. Thus, we have
\begin{equation*}\partial^{\alpha,\lambda}_{t} u(z,t)=e^{-\lambda t}\text{ }D_{t}^{\alpha} (e^{\lambda t} (u(z,t)-u(z,0)))-\lambda^{\alpha}(u(z,t)-u(z,0)),
\end{equation*}
where $\alpha\in(0,1)$, $\lambda>0$, and
$D_{t}^{\alpha}$ is the Caputo derivative defined as
\begin{equation*}
	D_{t}^{\alpha}g\left(t\right)=\frac{1}{\Gamma\left(1-\alpha\right)}\int_{0}^{t}\frac{d g\left(s\right)}{d s}\left(t-s\right)^{-\alpha}ds
\end{equation*}	for $g(t)\in AC[0,T]$ \cite{kilbas2006theory}.

\section{Finite difference method}
To solve the above problem numerically, we will approximate the limits by finite numbers and the derivatives by finite differences. After obtaining the discrete analogue of (\ref{12}) we will solve the problem recursively using boundary conditions.
\subsection{Weighted scheme for tsB-S model}
The system (\ref{12}) has the following form:
\begin{equation}
	\begin{cases}
		D_{t}^{\alpha}e^{\lambda t}(u\left(x,t\right)-u(x,0))= e^{\lambda t}(a\displaystyle\frac{\partial^2 u\left(x,t\right)}{\partial x^2}+b\displaystyle\frac{\partial u\left(x,t\right)}{\partial x}+(-c+\lambda^{\alpha})u\left(x,t\right)-u(x,0)\lambda^{\alpha} ),\\
		u\left(x,0\right)=f\left(x\right),\\
		u\left(x_{min},t\right)=p\left(t\right),\\
		u\left(x_{max},t\right)=q\left(t\right),\\
	\end{cases}\label{inwrs}
\end{equation}
where $a=\displaystyle\frac{1}{2}\sigma^2>0$, $b=\left(r-\displaystyle\frac{1}{2}\sigma^2\right)$, $c=r>0$, $\displaystyle f\left(x\right)=\max\left(e^{x}-K,0\right)$, $\displaystyle p\left(t\right)\rightarrow 0$ if $\displaystyle x_{min}\rightarrow -\infty$,
$\displaystyle q\left(t\right)\rightarrow e^{x_{max}}$ if $\displaystyle x_{max}\rightarrow \infty$. Since we assume that $x_{min}$ ($x_{max}$) is small (large) enough, we have $\displaystyle p\left(t\right)=0$. Moreover, the put-call parity (\ref{putcallpar}) implies that $\displaystyle q\left(t\right)=e^{x_{max}}-Ke^{-rt}$.  
Let us denote \[b_{j}=\left(j+1\right)^{1-\alpha}-j^{1-\alpha},\]
\[d=\Gamma\left(2-\alpha\right)\Delta t^{\alpha},\]
\[u^{k}=\left(u^{k}_{1},u^{k}_{2},\dots, u^{k}_{n-1}\right)^{T},\]
\[G^{k}=\left(\left(\displaystyle\frac{ad}{\Delta x^2}-\displaystyle\frac{bd}{2\Delta x}\right)u^{k}_{0},0,\dots ,0,\left(\displaystyle\frac{ad}{\Delta x^2}+\displaystyle\frac{bd}{2\Delta x}\right)u^{k}_{n}\right)^{T},\]
\[f=\left(f_{1},\dots ,f_{n-1}\right)^{T},\]
where $u^{k}_{i}=u\left(x_{i},t_{k}\right)$, $f_{i}=f\left(x_{i}\right)$, $i=1,\dots ,n-1$, $q^{k}=q\left(t_{k}\right)$, $k=0,1,\dots,N$, $j=0,1,\dots$, additionally $\displaystyle  \Delta t=T/N$, $\displaystyle  \Delta x=(x_{max}-x_{min})/n$ are time and space steps respectively. Moreover, $x_{i}=x_{min}+i\Delta x$ and $t_{j}=j\Delta t$ for $i=0,1\ldots,n$ and $j=0,1\ldots,N$, respectively. Note that $x_{min}$ and $x_{max}$ are the approximations of $-\infty$ and $\infty$, respectively. These values are determined experimentally such that further decreasing $x_{min}$ (or increasing $x_{max}$) has a negligible impact (e.g. of the machine epsilon's order) on the final solution. 
We will use the following approximations for space derivatives:
\begin{equation}
	\begin{cases}
		\displaystyle\frac{\partial u\left(x_{i},t_{k+1}\right)}{\partial x}=\displaystyle\frac{u\left(x_{i+1},t_{k+1}\right)-u\left(x_{i-1},t_{k+1}\right)}{2\Delta x}+O\left(\Delta x^{2}\right),\\
		\displaystyle\frac{\partial^{2} u\left(x_{i},t_{k+1}\right)}{\partial x^{2}}=\displaystyle\frac{u\left(x_{i+1},t_{k+1}\right)-2u\left(x_{i},t_{k+1}\right)+u\left(x_{i-1},t_{k+1}\right)}{\Delta x^{2}}+O\left(\Delta x^{2}\right),
	\end{cases}\label{space}
\end{equation}
where $k=0,1,\dots,N-1, l=1,\ldots,n-1$. We approximate the fractional-time derivative by \cite{Lin}:
\begin{equation}
	D_{t}^{\alpha}u\left(x_{i},t_{k+1}\right)=
	\displaystyle\frac{1}{\Gamma \left(2-\alpha \right)}\displaystyle\sum_{j=0}^{k}\displaystyle\frac{u\left(x_{i},t_{k+1-j}\right)-u\left(x_{i},t_{k-j}\right)}{\Delta t^{\alpha}}\left(\left(j+1\right)^{1-\alpha}-j^{1-\alpha}\right)+O\left(\Delta t^{2-\alpha}\right).\label{biolop}\end{equation}
After omitting the truncation errors, the implicit discrete scheme 
can be expressed in the following form: 

\begin{equation}
	\begin{cases}
		A\hat{u}^{1}=(1-d\lambda^{\alpha}) \hat{u}^{0}+G^{1},\\
		A\hat{u}^{k+1}=\displaystyle\sum_{j=0}^{k-1}\left(b_{j}-b_{j+1}\right)e^{-\left(j+1\right)\Delta t\lambda}(\hat{u}^{k-j}-\hat{u}^{0})+(1-d\lambda^{\alpha})\hat{u}^{0}+G^{k+1},
	\end{cases}\label{13}     
\end{equation}
where $k=1,\dots ,N-1$, $A=\left(a_{ij}\right)_{\left(n-1\right)\times\left(n-1\right)}$, such that:\newline
$a_{ij}=$
$$
\begin{cases}
	1+2\displaystyle\frac{ad}{\Delta x^2}+cd-\lambda^{\alpha}d,& \text{for } j=i,\text{ }i=1,2\dots ,n-1,\\
	-\left(\displaystyle\frac{ad}{\Delta x^2}-\displaystyle\frac{bd}{2\Delta x}\right),& \text{for } j=i-1,\text{ }i=2\dots ,n-1,\\
	-\left(\displaystyle\frac{ad}{\Delta x^2}+\displaystyle\frac{bd}{2\Delta x}\right),& \text{for } j=i+1,\text{ }i=1\dots ,n-2,\\
	0,& \text{in other cases}
\end{cases}
$$

The corresponding initial boundary conditions are as follows:
\begin{equation}
	\begin{cases}
		\hat{u}^{0}=f,\\
		\hat{u}^{k}_{0}=0,\\
		\hat{u}^{k}_{n}=q^{k},
	\end{cases}\label{14}     
\end{equation}
where $k=1,\dots ,N$. Similarly, let us write an explicit discrete scheme. We use approximations for space derivatives as follows:
\begin{equation}
	\begin{cases}
		\displaystyle\frac{\partial u\left(x_{i},t_{k}\right)}{\partial x}=\displaystyle\frac{u\left(x_{i+1},t_{k}\right)-u\left(x_{i-1},t_{k}\right)}{2\Delta x}+O\left(\Delta x^{2}\right),\\
		\displaystyle\frac{\partial^{2} u\left(x_{i},t_{k}\right)}{\partial x^{2}}=\displaystyle\frac{u\left(x_{i+1},t_{k}\right)-2u\left(x_{i},t_{k}\right)+u\left(x_{i-1},t_{k}\right)}{\Delta x^{2}}+O\left(\Delta x^{2}\right),
	\end{cases}\label{15}     
\end{equation}
where $k=0,1,\dots,N-1, l=1,\ldots,n-1$.

In matrix form, the explicit discrete scheme
can be expressed in the following form:
\begin{equation}
	\begin{cases}
		\hat{u}^{1}=e^{-\lambda \Delta t}B\hat{u}^{0}+(1-d\lambda^{\alpha}e^{-\lambda \Delta t}) \hat{u}^{0}+G^{0}e^{-\lambda \Delta t},\\
		\hat{u}^{k+1}=e^{-\lambda \Delta t}B\hat{u}^{k}+\displaystyle\sum_{j=0}^{k-1}\left(b_{j}-b_{j+1}\right)e^{-\left(j+1\right)\Delta t\lambda}(\hat{u}^{k-j}-\hat{u}^{0})+ (1-d\lambda^{\alpha}e^{-\lambda \Delta t}) \hat{u}^{0}+G^{k}e^{-\lambda \Delta t},
	\end{cases}\label{16}     
\end{equation}
where $k=1,\dots ,N-1$, $B=\left(b_{ij}\right)_{\left(n-1\right)\times\left(n-1\right)}$, such that:\newline
$b_{ij}=$
$
\begin{cases}
	-\left(2\displaystyle\frac{ad}{\Delta x^2}+cd-\lambda^{\alpha}d\right),& \text{for } j=i,\text{ }i=1\dots ,n-1,\\
	\displaystyle\frac{ad}{\Delta x^2}+\displaystyle\frac{bd}{2\Delta x},& \text{for } j=i+1,\text{ }i=1\dots ,n-2,\\
	\displaystyle\frac{ad}{\Delta x^2}-\displaystyle\frac{bd}{2\Delta x},& \text{for } j=i-1,\text{ }i=2\dots ,n-1,\\
	0,& \text{in other cases}
\end{cases}
$
\\Taking the linear combination of (\ref{13}) and (\ref{16}) we obtain a weighted scheme:

\begin{equation}
	\begin{cases}
		C\hat{u}^{1}=((1-\theta)(1-d\lambda^{\alpha})+\theta(1-d\lambda^{\alpha}e^{-\lambda \Delta t}))\hat{u}^{0}+\left(1-\theta\right)G^{1} +\theta G^{0}e^{-\lambda \Delta t}+\theta e^{-\lambda \Delta t} B\hat{u}^{0},\\
		C\hat{u}^{k+1}=\displaystyle\sum_{j=0}^{k-1}\left(b_{j}-b_{j+1}\right) e^{-\left(j+1\right)\Delta t\lambda}(\hat{u}^{k-j}-\hat{u}^{0})+\\((1-\theta)(1-d\lambda^{\alpha})+\theta(1-d\lambda^{\alpha}e^{-\lambda \Delta t}))\hat{u}^{0}+\left(1-\theta\right)G^{k+1} + \theta G^{k}e^{-\lambda \Delta t}+\theta e^{-\lambda \Delta t} B\hat{u}^{k},
	\end{cases}\label{17}       
\end{equation}
where $k=1,\dots ,N-1$, $C=\theta I+\left(1-\theta\right) A$, $\theta\in[0,1]$ 
and the corresponding initial boundary conditions are defined in (\ref{14}).


\subsection{Consistency of the weighted discrete scheme}
In this section, we show the following 
\begin{Theorem}
	For $\theta\in[0,1]$ and $1\leq i\leq n-1,1\leq j\leq N$, the truncation error $R_{i}^{j}$ of the numerical scheme (\ref{17}) satisfies
	\[\left|R_{i}^{j}\right|\leq C_{1}\Delta t^{\alpha}\left(\Delta t^{2-\alpha}+\Delta x^{2}\right).\]
	\label{rukuio}\end{Theorem}

\begin{Proof}
	Let us replace the derivatives in the first equation of (\ref{inwrs}) by their discrete equivalents. Then, by (\ref{space}), (\ref{biolop}) and (\ref{15}), we get the full formulation of the discrete weighted scheme with the truncation error:
	\begin{equation}
		\begin{cases}
			-\left(\displaystyle\frac{ad}{\Delta x^{2}}+\displaystyle\frac{bd}{2\Delta x}\right)(\theta u^{0}_{i+1}+(1-\theta) u^{1}_{i+1})+\left(\displaystyle\frac{2ad}{\Delta x^{2}}+cd-d\lambda^{\alpha}\right)(\theta u^{0}_{i}+(1-\theta) u^{1}_{i})-\\ \left(\displaystyle\frac{ad}{\Delta x^{2}}-\displaystyle\frac{bd}{2\Delta x}\right)(\theta u^{0}_{i-1}+(1-\theta) u^{1}_{i-1})=  (\theta e^{\lambda \Delta t}+(1-\theta)-d\lambda^{\alpha})u_{i}^{0}-(\theta e^{\lambda \Delta t}+1-\theta)u_{i}^{1}+(1-\theta)R_{i}^{1}+\theta R_{i}^{0},\\
			-\left(\displaystyle\frac{ad}{\Delta x^{2}}+\displaystyle\frac{bd}{2\Delta x}\right)(\theta u^{k}_{i+1}+(1-\theta) u^{k+1}_{i+1})+\left(\displaystyle\frac{2ad}{\Delta x^{2}}+cd-d\lambda^{\alpha}\right)(\theta u^{k}_{i}+(1-\theta) u^{k+1}_{i})-\\  \left(\displaystyle\frac{ad}{\Delta x^{2}}-\displaystyle\frac{bd}{2\Delta x}\right)(\theta u^{k}_{i-1}+(1-\theta) u^{k+1}_{i-1})=
			\displaystyle\sum_{j=0}^{k-1}\left(b_{j}-b_{j+1}\right)((1-\theta) e^{-\left(j+1\right)\Delta t\lambda}+\theta e^{-j\Delta t\lambda})(u_{i}^{k-j}-u_{i}^{0})+\\(\theta e^{\lambda \Delta t}+(1-\theta)-d\lambda^{\alpha})u_{i}^{0}-
			(\theta e^{\lambda \Delta t}+1-\theta)u^{k+1}_{i}+(1-\theta)R_{i}^{k+1}+\theta R_{i}^{k} .
		\end{cases}\label{18}     
	\end{equation}
	Here $k=1,\dots N-1$, $i=1,\dots ,n-1$, $R_{i}^{k}=d\varrho_{i}^{k}$ for the (original) truncation error $\varrho_{i}^{k}$ and the corresponding initial boundary conditions are defined in (\ref{14}). Note that in this paper, both $R_{i}^{k}$ and $\varrho_{i}^{k}$ are called truncation errors, but only $R_{i}^{k}$ will be considered further. 
	By (\ref{space}), (\ref{biolop}) and (\ref{15}), we have
	\[\left|R_{i}^{j}\right|\leq C_{i}^{j}\Delta t^{\alpha}(\Delta t^{2-\alpha}+\Delta x^{2}),\]
	where $C_{i}^{j}$ are constants ($1\leq i\leq n,1\leq j\leq N$). Let us denote $\displaystyle C_{1}=\max_{1\leq i\leq n,1\leq j\leq N}C_{i}^{j}$.
	Then, for the truncation error, it holds that
	\[\left|R_{i}^{j}\right|\leq C_{1}\Delta t^{\alpha}\left(\Delta t^{2-\alpha}+\Delta x^{2}\right).\]
\end{Proof}

Note that the parameters $\theta$ and $\lambda$ have no influence in the above analysis.

\subsection{Stability of the weighted discrete scheme}
For the purposes of stability and convergence analysis, for $l=0,1,\dots n$ and $k=0,1,..,N$ let us denote:
$u_{l}^{k}=u(x_{l},t_{k})$ - the exact solution of (\ref{12}) evaluated at the grid point,
$\hat{u}_{l}^{k}=\hat{u}\left(x_{l},t_{k}\right)$ - the exact solution of the numerical scheme (\ref{17}),
$\tilde{u}_{l}^{k}$ - some approximation of $\hat{u}_{l}^{k}$. We will proceed using the von Neumann method.
After omitting the truncation error and introducing $\epsilon_{l}^{k}=\hat{u}_{l}^{k}-\tilde{u}_{l}^{k}$ ($l=0,1,\dots n$, $k=0,1,..,N$), (\ref{18}) can be transformed into: 
\begin{equation}
	\begin{cases}
		-\left(\displaystyle\frac{ad}{\Delta x^{2}}+\displaystyle\frac{bd}{2\Delta x}\right)(\theta \epsilon^{0}_{i+1}+(1-\theta) \epsilon^{1}_{i+1})+\left(\displaystyle\frac{2ad}{\Delta x^{2}}+cd-d\lambda^{\alpha}\right)(\theta \epsilon^{0}_{i}+(1-\theta) \epsilon^{1}_{i})- \left(\displaystyle\frac{ad}{\Delta x^{2}}-\displaystyle\frac{bd}{2\Delta x}\right)(\theta \epsilon^{0}_{i-1}+(1-\theta) \epsilon^{1}_{i-1})=\\  (\theta \epsilon^{\lambda \Delta t}+(1-\theta)-d\lambda^{\alpha})\epsilon_{i}^{0}-(\theta e^{\lambda \Delta t}+1-\theta)\epsilon_{i}^{1},\\
		-\left(\displaystyle\frac{ad}{\Delta x^{2}}+\displaystyle\frac{bd}{2\Delta x}\right)(\theta \epsilon^{k}_{i+1}+(1-\theta) \epsilon^{k+1}_{i+1})+\left(\displaystyle\frac{2ad}{\Delta x^{2}}+cd-d\lambda^{\alpha}\right)(\theta \epsilon^{k}_{i}+(1-\theta) \epsilon^{k+1}_{i})-  \left(\displaystyle\frac{ad}{\Delta x^{2}}-\displaystyle\frac{bd}{2\Delta x}\right)(\theta \epsilon^{k}_{i-1}+(1-\theta) \epsilon^{k+1}_{i-1})=\\
		\displaystyle\sum_{j=0}^{k-1}\left(b_{j}-b_{j+1}\right)(\theta e^{-j\Delta t\lambda}+(1-\theta) e^{-\left(j+1\right)\Delta t\lambda})(\epsilon_{i}^{k-j}-\epsilon_{i}^{0})+(\theta e^{\lambda \Delta t}+(1-\theta)-d\lambda^{\alpha})\epsilon_{i}^{0}-\\(\theta e^{\lambda \Delta t}+1-\theta)\epsilon^{k+1}_{i},\\
		\epsilon_{0}^{m}=\epsilon_{n}^{m}=0,
	\end{cases}\label{19}     
\end{equation}
where $k=1,\dots ,N-1$, $m=1,\dots ,N$, $i=1,\dots ,n-1$. Note that the last equality in (\ref{19}) holds because we assume that $\hat{u}_{0}^{m}=\tilde{u}_{0}^{m}$ and $\hat{u}_{n}^{m}=\tilde{u}_{n}^{m}$.
Now, let us introduce the following grid function:\\
$\displaystyle \epsilon^{k}(x)=$
\begin{equation*}
	\begin{cases}
		\displaystyle \epsilon_{l}^{k}\text{, }x\in\left(x_{l}-\Delta x/2,x_{l}+\Delta x/2\right],\text{ }l=1,2\dots ,n-1,\\
		\displaystyle 0\text{, }x\in\left(x_{min},x_{min}+\Delta x/2\right]\cup\left[x_{max}-\Delta x/2,x_{max}\right].
	\end{cases}
\end{equation*}
Because $\displaystyle \epsilon^{k}_{0}=\epsilon^{k}_{n}$, we make a periodic expansion for $\epsilon^{k}_{l}$ with period $Y=x_{max}-x_{min}$. 
Then $\displaystyle \epsilon^{k}(x)$ has the following Fourier series extension:
\begin{equation*}
	\displaystyle \epsilon^{k}(x)=\displaystyle\sum_{j=-\infty}^{\infty}v_{j}^{k}e^{2j\pi x \iota/Y},
\end{equation*}
where $\displaystyle v_{j}^{k}=\frac{1}{Y}\int_{0}^{Y}\epsilon^{k}(x)e^{2j\pi x \iota/Y}dx$, $\iota=\sqrt{-1}$, $k=0, 1,\dots N$.
We define the norm $\norm{\cdot}_{\Delta x}$ as \[\norm{\epsilon^{k}}_{\Delta x}=\sqrt{\sum_{j=1}^{n-1}\Delta x\left|\epsilon_{j}^{k}\right|^{2}},\] where 
\[\epsilon^{k}=(\epsilon_{1}^{k},\epsilon^{k}_{2},\dots ,\epsilon_{n-1}^{k}). \]
Because $\epsilon_{0}^{k}=\epsilon_{n}^{k}=0,$ it follows
\[\norm{\epsilon^{k}}^{2}_{\Delta x}=\int_{0}^{Y}\left|\epsilon^{k}(x)\right|^{2}dx=\norm{\epsilon^{k}(x)}^{2},\]
where $\norm{\cdot}$ is $L^{2}[0,Y].$ Using the Parseval identity we have:
\[\norm{\epsilon^{k}}^{2}_{\Delta x}=\sum_{j=1}^{n-1}\Delta x\left|\epsilon_{j}^{k}\right|^{2}=Y\sum_{j=-\infty}^{\infty}\left|v_{j}^{k}\right|^{2},\]
$k=0,1,\dots ,N.$ Based on the above analysis and the fact that $x_{l}=x_{min}+l\Delta x$ for $l=0,\ldots,n$, we 
infer that the solution of (\ref{19}),
has the form:
\begin{equation}
	\displaystyle \epsilon_{l}^{k}=v^{k}e^{\iota\beta (x_{min}+l\Delta x)}, \label{20}     
\end{equation}
where $\displaystyle \beta=\frac{2\pi j}{Y}.$
Substituting into (\ref{19}) we get:
\begin{equation}
	\begin{cases}
		\displaystyle \left(-\left(\frac{ad}{\Delta x^{2}}+\frac{bd}{2\Delta x}\right)e^{\iota\beta \Delta x}+\left(2\frac{ad}{\Delta x^{2}}+cd-d\lambda^{\alpha}\right)-\left(\frac{ad}{\Delta x^{2}}-\frac{bd}{2\Delta x}\right)
		e^{-\iota\beta \Delta x}\right)\left(\theta v^{0}+\left(1-\theta\right)v^{1}\right)=\\(\theta e^{\lambda \Delta t}+(1-\theta)-d\lambda^{\alpha})v^{0}-(\theta e^{\lambda \Delta t}+1-\theta)v^{1},\\
		\displaystyle \left(-\left(\frac{ad}{\Delta x^{2}}
		+\frac{bd}{2\Delta x}\right)e^{\iota\beta \Delta x}+
		\left(2\frac{ad}{\Delta x^{2}}+cd-d\lambda^{\alpha}\right)
		-\left(\frac{ad}{\Delta x^{2}}-\frac{bd}{2\Delta x}\right)
		e^{-\iota\beta \Delta x}\right)
		\displaystyle\left(\theta v^{k}+\left(1-\theta\right)v^{k+1}\right)=\\
		\displaystyle\sum_{j=0}^{k-1}\left(b_{j}-b_{j+1}\right)(\theta e^{-j\Delta t\lambda}+(1-\theta) e^{-\left(j+1\right)\Delta t\lambda})(v^{k-j}-v^{0})+(\theta e^{\lambda \Delta t}+(1-\theta)-d\lambda^{\alpha})v^{0}-\\(\theta e^{\lambda \Delta t}+1-\theta)v^{k+1},\\
	\end{cases}\label{21}     
\end{equation}
for $k=1,\dots ,N-1$. To continue, we have to find a relation between the coefficients $b_{j}$.

\begin{Proposition}
	For $N\geq1$ we have
	\begin{equation}
		\frac{1-\alpha}{(N+1)^{\alpha}}<(N+1)^{1-\alpha}-N^{1-\alpha}<\frac{1-\alpha}{N^{\alpha}}.\label{1-alf}
	\end{equation}\label{jmp50krol}
\end{Proposition}
\begin{Proof}
	Let us introduce the function $f(x)=x^{1-\alpha}$ for $x>0$. By Lagrange's Mean Value Theorem, we have
	\[\frac{f(x+1)-f(x)}{(x+1)-x}=f'(C)=\frac{1-\alpha}{C^{\alpha}},\]
	where $C\in(x,x+1)$. Since $f'$ is a decreasing function, the proof is completed.
\end{Proof}
\begin{Proposition}\cite{ja2}\label{3.1}
	The coefficients $b_{j}=(j+1)^{1-\alpha}-j^{1-\alpha}$ satisfy:
	\begin{enumerate}
		\item $b_{j}>0,$ $j=0,1\dots $
		\item $1=b_{0}>b_{1}>\dots >b_{k}$
		\item $\displaystyle\lim_{k\to\infty}b_{k}=0$
		\item $\displaystyle\sum_{j=0}^{k-1}(b_{j}-b_{j+1})+b_{k}=1$
	\end{enumerate}
\end{Proposition} 
\begin{Proof}
	\begin{enumerate}
		\item $b_{j}=(j+1)^{1-\alpha}-j^{1-\alpha}>0,$ for $j\geq0$ and $\alpha\in(0,1)$. 
		\item This is the conclusion of Proposition \ref{jmp50krol}.
		\item This is the conclusion of Proposition \ref{jmp50krol}.
		\item $\displaystyle\sum_{j=0}^{k-1}(b_{j}-b_{j+1})+b_{k}=(1-b_{1})+(b_{1}-b_{2})+(b_{2}-b_{3})+\dots+(b_{k-1}-b_{k})+b_{k}=1.$
	\end{enumerate}
\end{Proof}

Now we will check under which conditions $\left|v^{m}\right|\leq\left|v^{0}\right|$ for each $m=1,\dots ,N$. Then $\norm{\epsilon^{m}}\leq \norm{\epsilon^{0}}$,
in other words, the weighted scheme is stable.
\begin{Theorem}
	Let $\theta\in[0,1]$. If the following inequality
	
	\begin{dmath}
		\displaystyle\min\left( \displaystyle\left|\left(1-\theta\right)\left(cd-d\lambda^{\alpha}+1\right)+\theta e^{\lambda \Delta t}\right|,\left|\left(1-\theta\right)\left(\frac{4ad}{\Delta x^{2}}+cd-d\lambda^{\alpha}+1\right)+\theta e^{\lambda \Delta t}\right|\right)\geq
		\sqrt{\max\left(\left(\left(1-b_{1}\right)\left(\theta+\left(1-\theta\right)e^{-\lambda\Delta t}\right)-\frac{4ad\theta}{\Delta x^{2}}-cd\theta +\lambda^{\alpha}d\theta\right)^{2}, \left(\left(1-b_{1}\right)\left(\theta+\left(1-\theta\right)e^{-\lambda\Delta t}\right)-cd\theta +\lambda^{\alpha}d\theta\right)^{2}+\left(\frac{bd\theta}{\Delta x}\right)^{2}\right)}+\sum_{j=1}^{N-1}\left(b_{j}-b_{j+1}\right)\left(\theta e^{-\lambda j \Delta t}+\left(1-\theta\right) e^{-\lambda \left(j+1\right) \Delta t}\right)+\left|\theta e^{\lambda\Delta t}+1-\theta-d\lambda^{\alpha}-\sum_{j=0}^{N-1}\left(b_{j}-b_{j+1}\right)\left(\theta e^{-\lambda j \Delta t}+\left(1-\theta\right) e^{-\lambda \left(j+1\right) \Delta t}\right)\right|\label{22}  
	\end{dmath}   
	
	holds, then the scheme (\ref{17}) is stable.
	
\end{Theorem}

\begin{Proof}
	We have to show that $v^{m}$ defined in (\ref{20}) follows
	$\displaystyle\left|v^{m}\right|\leq\left|v^{0}\right|$ for $m=2,3,\dots ,k.$ 
	Let us denote\\
	\begin{multline*}
		\zeta=\left(-4\sin^{2}\left(\frac{\beta \Delta x}{2}\right)+2\right)\left(\frac{-ad}{\Delta x^{2}}\right)+2\frac{ad}{\Delta x^{2}}+cd-2\iota\frac{bd}{2\Delta x}\sin\left(\beta \Delta x\right)=
		\sin^{2}\left(\frac{\beta \Delta x}{2}\right)\frac{4ad}{\Delta x^{2}}+cd-\iota\frac{bd}{\Delta x}\sin\left(\beta \Delta x\right).\end{multline*}
	Let us observe that $\displaystyle Re{\text{ }\zeta}=\sin^{2}\left(\frac{\beta \Delta x}{2}\right)\frac{4ad}{\Delta x^{2}}+cd>0$. The proof of this fact is immediate because $a,d,c,\Delta x>0$.
	
	At the beginning, we will show that $\forall 0<k<N$ we have
	\begin{multline}
		\sum_{j=1}^{k}\left(b_{j}-b_{j+1}\right)\left(\theta e^{-\lambda j \Delta t}+\left(1-\theta\right) e^{-\lambda \left(j+1\right) \Delta t}\right)+
		\left|\theta e^{\lambda \Delta t}+1-\theta-d\lambda^{\alpha}-\sum_{j=0}^{k}\left(b_{j}-b_{j+1}\right)\left(\theta e^{-\lambda j \Delta t}+\left(1-\theta\right) e^{-\lambda \left(j+1\right) \Delta t}\right)\right|\leq\\
		\sum_{j=1}^{N-1}\left(b_{j}-b_{j+1}\right)\left(\theta e^{-\lambda j \Delta t}+\left(1-\theta\right) e^{-\lambda \left(j+1\right) \Delta t}\right)+
		\left|\theta e^{\lambda \Delta t}+1-\theta-d\lambda^{\alpha}-\sum_{j=0}^{N-1}\left(b_{j}-b_{j+1}\right)\left(\theta e^{-\lambda j \Delta t}+\left(1-\theta\right) e^{-\lambda \left(j+1\right) \Delta t}\right)\right|.\label{23}     
	\end{multline}
	
	It is equivalent to
	
	\begin{multline}
		0\leq
		\sum_{j=k+1}^{N-1}\left(b_{j}-b_{j+1}\right)\left(\theta e^{-\lambda j \Delta t}+\left(1-\theta\right) e^{-\lambda \left(j+1\right) \Delta t}\right)+\\
		\left|\theta e^{\lambda \Delta t}+1-\theta-d\lambda^{\alpha}-\sum_{j=0}^{N-1}\left(b_{j}-b_{j+1}\right)\left(\theta e^{-\lambda j \Delta t}+\left(1-\theta\right) e^{-\lambda \left(j+1\right) \Delta t}\right)\right|-\\
		\left|\theta e^{\lambda \Delta t}+1-\theta-d\lambda^{\alpha}-\sum_{j=0}^{k}\left(b_{j}-b_{j+1}\right)\left(\theta e^{-\lambda j \Delta t}+\left(1-\theta\right) e^{-\lambda \left(j+1\right) \Delta t}\right)\right|.\label{24}     
	\end{multline}
	Let us consider $3$ cases:
	\begin{enumerate}
		\item  $\displaystyle\theta e^{\lambda \Delta t}+1-\theta-d\lambda^{\alpha}\geq \sum_{j=0}^{N-1}\left(b_{j}-b_{j+1}\right)\left(\theta e^{-\lambda j \Delta t}+\left(1-\theta\right) e^{-\lambda \left(j+1\right) \Delta t}\right).$ 
		Then (\ref{24}) has the form: 
		\begin{multline*}
			0\leq \sum_{j=k+1}^{N-1}\left(b_{j}-b_{j+1}\right)\left(\theta e^{-\lambda j \Delta t}+\left(1-\theta\right) e^{-\lambda \left(j+1\right) \Delta t}\right)+\sum_{j=0}^{k}\left(b_{j}-b_{j+1}\right)\left(\theta e^{-\lambda j \Delta t}+\left(1-\theta\right) e^{-\lambda \left(j+1\right) \Delta t}\right)-\\ \sum_{j=0}^{N-1}\left(b_{j}-b_{j+1}\right)\left(\theta e^{-\lambda j \Delta t}+\left(1-\theta\right) e^{-\lambda \left(j+1\right) \Delta t}\right)=0.
		\end{multline*}
		
		\item  \begin{equation}
			\begin{cases}
				\displaystyle\theta e^{\lambda \Delta t}+1-\theta-d\lambda^{\alpha}\geq \sum_{j=0}^{k}\left(b_{j}-b_{j+1}\right)\left(\theta e^{-\lambda j \Delta t}+\left(1-\theta\right) e^{-\lambda \left(j+1\right) \Delta t}\right),\\
				\displaystyle\theta e^{\lambda \Delta t}+1-\theta-d\lambda^{\alpha}\leq \sum_{j=0}^{N-1}\left(b_{j}-b_{j+1}\right)\left(\theta e^{-\lambda j \Delta t}+\left(1-\theta\right) e^{-\lambda \left(j+1\right) \Delta t}\right).
			\end{cases}\label{25}     
		\end{equation}
		Then (\ref{24}) has the form: 
		\begin{multline*}
			0\leq \sum_{j=k+1}^{N-1}\left(b_{j}-b_{j+1}\right)\left(\theta e^{-\lambda j \Delta t}+\left(1-\theta\right) e^{-\lambda \left(j+1\right) \Delta t}\right)+\sum_{j=0}^{k}\left(b_{j}-b_{j+1}\right)\left(\theta e^{-\lambda j \Delta t}+\left(1-\theta\right) e^{-\lambda \left(j+1\right) \Delta t}\right)+\\ \sum_{j=0}^{N-1}\left(b_{j}-b_{j+1}\right)\left(\theta e^{-\lambda j \Delta t}+\left(1-\theta\right) e^{-\lambda \left(j+1\right) \Delta t}\right)-2\left(\theta e^{\lambda \Delta t}+1-\theta-d\lambda^{\alpha}\right)=\\
			2\left(\sum_{j=0}^{N-1}\left(b_{j}-b_{j+1}\right)\left(\theta e^{-\lambda j \Delta t}+\left(1-\theta\right) e^{-\lambda \left(j+1\right) \Delta t}\right)-\left(\theta e^{\lambda \Delta t}+1-\theta-d\lambda^{\alpha}\right)\right),
		\end{multline*}
		where the expression on the right side is always positive by the first equation of (\ref{25}).
		
		\item  $\displaystyle\theta e^{\lambda \Delta t}+1-\theta-d\lambda^{\alpha}\leq \sum_{j=0}^{k}\left(b_{j}-b_{j+1}\right)\left(\theta e^{-\lambda j \Delta t}+\left(1-\theta\right) e^{-\lambda \left(j+1\right) \Delta t}\right).$ 
		Then (\ref{24}) has the form: 
		\begin{multline*}
			0\leq \sum_{j=k+1}^{N-1}\left(b_{j}-b_{j+1}\right)\left(\theta e^{-\lambda j \Delta t}+\left(1-\theta\right) e^{-\lambda \left(j+1\right) \Delta t}\right)+\sum_{j=0}^{N-1}\left(b_{j}-b_{j+1}\right)\left(\theta e^{-\lambda j \Delta t}+\left(1-\theta\right) e^{-\lambda \left(j+1\right) \Delta t}\right)-\\\left(\theta e^{\lambda \Delta t}+1-\theta-d\lambda^{\alpha}\right)+\left(\theta e^{\lambda \Delta t}+1-\theta-d\lambda^{\alpha}\right)-\sum_{j=0}^{k}\left(b_{j}-b_{j+1}\right)\left(\theta e^{-\lambda j \Delta t}+\left(1-\theta\right) e^{-\lambda \left(j+1\right) \Delta t}\right)=\\
			2\sum_{j=k+1}^{N-1}\left(b_{j}-b_{j+1}\right)\left(\theta e^{-\lambda j \Delta t}+\left(1-\theta\right) e^{-\lambda \left(j+1\right) \Delta t}\right),
		\end{multline*} 
		where the expression on the right side is always positive by Proposition \ref{3.1}. In this way, we have shown that (\ref{23}) holds. 
	\end{enumerate} 
	
	Now we will follow the mathematical induction method to show that for each $m=1,2\dots,N$ we have $\left|v^{m}\right|\leq\left|v^{0}\right|$.
	\begin{enumerate}
		\item $n=1$
		By the identity \[ \sin^{2}\frac{z}{2}=-\frac{1}{4}\left(e^{\iota z}-2+e^{-\iota z}\right),\]
		the first equation of (\ref{21}) can be transformed into
		\[\zeta\left(\left(1-\theta\right)v^{1}+\theta v^{0}\right)-d\lambda^{\alpha}\left(\left(1-\theta\right)v^{1}+\theta v^{0}\right)=v^{0}(\theta e^{\lambda\Delta t}+1-\theta-d\lambda^{\alpha})-v^{1}(\theta e^{\lambda\Delta t}+1-\theta),\]
		
		So 
		\[\displaystyle\left|\zeta\left(1-\theta\right)-(1-\theta)d\lambda^{\alpha}+\theta e^{\lambda \Delta t}+1-\theta\right|\left|v^{1}\right|=\left|\theta e^{\lambda \Delta t}+1-\theta-d\lambda^{\alpha}-\zeta\theta+d\lambda^{\alpha}\theta \right| \left|v^{0}\right|.\]
		If
		\begin{equation*}
			\displaystyle\left|\theta e^{\lambda \Delta t}+1-\theta-d\lambda^{\alpha}-\zeta\theta+d\lambda^{\alpha}\theta \right|\leq \left|\zeta\left(1-\theta\right)-(1-\theta)d\lambda^{\alpha}+\theta e^{\lambda \Delta t}+1-\theta\right|.
		\end{equation*} 
		then \[\left|v^{1}\right|\leq\left|v^{0}\right|.\]
		Note, that
		\begin{multline*}
			\displaystyle\left|\theta e^{\lambda \Delta t}+1-\theta-d\lambda^{\alpha}-\zeta\theta+d\lambda^{\alpha}\theta \right|=\\
			\left|-\zeta\theta+d\lambda^{\alpha}\theta+\left(1-b_{1}\right)\left(\theta+\left(1-\theta\right)e^{-\lambda\Delta t}\right)+\theta e^{\lambda \Delta t}+1-\theta-d\lambda^{\alpha}-\left(1-b_{1}\right)\left(\theta+\left(1-\theta\right)e^{-\lambda\Delta t}\right)\right|
			\leq\\
			\left|-\zeta\theta+d\lambda^{\alpha}\theta +\left(1-b_{1}\right)\left(\theta+\left(1-\theta\right)e^{-\lambda\Delta t}\right)\right|+\left| \theta e^{\lambda \Delta t}+1-\theta-d\lambda^{\alpha}-\left(1-b_{1}\right)\left(\theta+\left(1-\theta\right)e^{-\lambda\Delta t}\right) \right|\leq\\         
			\sqrt{(-\theta(\sin^{2}\left(\frac{\beta \Delta x}{2}\right)\frac{4ad}{\Delta x^{2}}+cd)+ d\lambda^{\alpha}\theta+\left(1-b_{1}\right)\left(\theta+\left(1-\theta\right)e^{-\lambda\Delta t}\right))^{2}+
				(\theta\frac{bd}{\Delta x}\sin(\beta \Delta x))^{2}}+\\
			\left|\sum_{j=1}^{k-1}\left(b_{j}-b_{j+1}\right)\left(\theta e^{-\lambda j \Delta t}+\left(1-\theta\right) e^{-\lambda \left(j+1\right) \Delta t}\right)+\theta e^{\lambda \Delta t}+1-\theta-d\lambda^{\alpha}-\sum_{j=0}^{k-1}\left(b_{j}-b_{j+1}\right)\left(\theta e^{-\lambda j \Delta t}+\left(1-\theta\right) e^{-\lambda \left(j+1\right) \Delta t}\right)\right|\leq\\
			\sqrt{\max\left(\left(\frac{4ad\theta}{\Delta x^{2}}+cd\theta -\lambda^{\alpha}d\theta-\left(1-b_{1}\right)\left(\theta+\left(1-\theta\right)e^{-\lambda\Delta t}\right)\right)^{2}, \left(cd\theta -\lambda^{\alpha}d\theta-\left(1-b_{1}\right)\left(\theta+\left(1-\theta\right)e^{-\lambda\Delta t}\right)\right)^{2}\right)+\left(\frac{bd\theta}{\Delta x}\right)^{2}}+\\
			\sum_{j=1}^{k-1}\left(b_{j}-b_{j+1}\right)\left(\theta e^{-\lambda j \Delta t}+\left(1-\theta\right) e^{-\lambda \left(j+1\right) \Delta t}\right)+
			\left|\theta e^{\lambda \Delta t}+1-\theta-d\lambda^{\alpha}-\sum_{j=0}^{k-1}\left(b_{j}-b_{j+1}\right)\left(\theta e^{-\lambda j \Delta t}+\left(1-\theta\right) e^{-\lambda \left(j+1\right) \Delta t}\right)\right|\leq\\
			\sqrt{\max\left(\left(\frac{4ad\theta}{\Delta x^{2}}+cd\theta -\lambda^{\alpha}d\theta-\left(1-b_{1}\right)\left(\theta+\left(1-\theta\right)e^{-\lambda\Delta t}\right)\right)^{2}, \left(cd\theta -\lambda^{\alpha}d\theta-\left(1-b_{1}\right)\left(\theta+\left(1-\theta\right)e^{-\lambda\Delta t}\right)\right)^{2}\right)+\left(\frac{bd\theta}{\Delta x}\right)^{2}}+\\
			\sum_{j=1}^{N-1}\left(b_{j}-b_{j+1}\right)\left(\theta e^{-\lambda j \Delta t}+\left(1-\theta\right) e^{-\lambda \left(j+1\right) \Delta t}\right)+
			\left|\theta e^{\lambda \Delta t}+1-\theta-d\lambda^{\alpha}-\sum_{j=0}^{N-1}\left(b_{j}-b_{j+1}\right)\left(\theta e^{-\lambda j \Delta t}+\left(1-\theta\right) e^{-\lambda \left(j+1\right) \Delta t}\right)\right|,         \end{multline*} 
		where the last inequality is true by (\ref{23}).
		On the other hand, 
		\begin{multline}
			\left|\left(1-\theta\right)\left(\zeta-d\lambda^{\alpha}\right)+\theta e^{\lambda \Delta t}+1-\theta\right|=\sqrt{\left(\left(1-\theta\right)\left(\sin^{2}\frac{\beta\Delta x}{2}\frac{4ad}{\Delta x^{2}}+cd-d\lambda^{\alpha}+1\right)+\theta e^{\lambda\Delta t}\right)^{2}+\left(\left(1-\theta\right)\frac{bd}{\Delta x}\sin\beta\Delta x\right)^{2}}\geq\\
			\displaystyle\sqrt{\min\left(\left(\left(1-\theta\right)\left(cd-d\lambda^{\alpha}+1\right)+\theta e^{\lambda \Delta t}\right)^{2},\left(\left(1-\theta\right)\left(\frac{4ad}{\Delta x^{2}}+cd-d\lambda^{\alpha}+1\right)+\theta e^{\lambda \Delta t}\right)^{2}\right)}=\\
			\min\left( \displaystyle\left|\left(1-\theta\right)\left(cd-d\lambda^{\alpha}+1\right)+\theta e^{\lambda \Delta t}\right|,\left|\left(1-\theta\right)\left(\frac{4ad}{\Delta x^{2}}+cd-d\lambda^{\alpha}+1\right)+\theta e^{\lambda \Delta t}\right|\right).\label{26}     
		\end{multline}
		So, as the conlusion, (\ref{22}) implies 
		\begin{equation*}
			\displaystyle\left|\theta e^{\lambda \Delta t}+1-\theta-d\lambda^{\alpha}-\zeta\theta+d\lambda^{\alpha}\theta \right|\leq \left|\zeta\left(1-\theta\right)-(1-\theta)d\lambda^{\alpha}+\theta e^{\lambda \Delta t}+1-\theta\right|,
		\end{equation*} 
		from which it follows: \[\left|v^{1}\right|\leq\left|v^{0}\right|.\]
		
		\item Let us suppose that
		\[\displaystyle\left|v^{m}\right|\leq\left|v^{0}\right|,\]
		for $m=1,2,\dots ,k,$ $k<N.$\\
		To complete the proof, we have to show that 
		\[\displaystyle\left|v^{k+1}\right|\leq\left|v^{0}\right|.\]
		By the second equation of (\ref{21})
		we get 
		\begin{multline*}
			|v^{k+1}\left(\zeta\left(1-\theta\right)-\left(1-\theta\right)d\lambda^{\alpha}+\theta e^{\lambda \Delta t}+1-\theta\right)|=\\|v^{k}\left(-\theta\zeta +d\lambda^{\alpha}\theta\right)+\sum_{j=0}^{k-1}\left(b_{j}-b_{j+1}\right)\left(\theta e^{-\lambda j \Delta t}+\left(1-\theta\right) e^{-\lambda \left(j+1\right) \Delta t}\right)\left(v^{k-j}-v^{0}\right)+\left(\theta e^{\lambda \Delta t}+1-\theta-d\lambda^{\alpha}\right)v^{0}|.
		\end{multline*}
		
		Note that
		\begin{equation} 
			|v^{k+1}\left(\zeta\left(1-\theta\right)-\left(1-\theta\right)d\lambda^{\alpha}+\theta e^{\lambda \Delta t}+1-\theta\right)|=
			|v^{k+1}||\left(\zeta\left(1-\theta\right)-\left(1-\theta\right)d\lambda^{\alpha}+\theta e^{\lambda \Delta t}+1-\theta\right)|,\label{27}     
		\end{equation} 
		and\\
		\begin{multline}\displaystyle\left|v^{k}\left(-\theta\zeta +d\lambda^{\alpha}\theta\right)+\sum_{j=0}^{k-1}\left(b_{j}-b_{j+1}\right)\left(\theta e^{-\lambda j \Delta t}+\left(1-\theta\right) e^{-\lambda \left(j+1\right) \Delta t}\right)\left(v^{k-j}-v^{0}\right)+\left(\theta e^{\lambda \Delta t}+1-\theta-d\lambda^{\alpha}\right)v^{0}\right|
			\leq\\
			\left|v^{k}\right| \left|-\theta\zeta +d\lambda^{\alpha}\theta+\left(1-b_{1}\right)\left(\theta+\left(1-\theta\right)e^{-\lambda\Delta t}\right)\right| 
			+\left|\sum_{j=1}^{k-1}\left(b_{j}-b_{j+1}\right)\left(\theta e^{-\lambda j \Delta t}+
			\left(1-\theta\right) e^{-\lambda \left(j+1\right) \Delta t}\right)\right|
			\left|v^{k-j}\right|+\\
			\left|\theta e^{\lambda \Delta t}+1-\theta-d\lambda^{\alpha}-\sum_{j=0}^{k-1}\left(b_{j}-b_{j+1}\right)\left(\theta e^{-\lambda j \Delta t}+\left(1-\theta\right) e^{-\lambda \left(j+1\right) \Delta t}\right)\right|\left|v^{0}\right|
			\leq\\
			\left|v^{0}\right|\sqrt{(-\theta(\sin^{2}\left(\frac{\beta \Delta x}{2}\right)\frac{4ad}{\Delta x^{2}}+cd)+ d\lambda^{\alpha}\theta+\left(1-b_{1}\right)\left(\theta+\left(1-\theta\right)e^{-\lambda\Delta t}\right))^{2}+
				(\theta\frac{bd}{\Delta x}\sin(\beta \Delta x))^{2}}+\\
			\left|v^{0}\right| \sum_{j=1}^{k-1}\left(b_{j}-b_{j+1}\right)\left(\theta e^{-\lambda j \Delta t}+\left(1-\theta\right) e^{-\lambda \left(j+1\right) \Delta t}\right)+\left|v^{0}\right|\left|\theta e^{\lambda \Delta t}+1-\theta-d\lambda^{\alpha}-\sum_{j=0}^{k-1}\left(b_{j}-b_{j+1}\right)\left(\theta e^{-\lambda j \Delta t}+\left(1-\theta\right) e^{-\lambda \left(j+1\right) \Delta t}\right)\right|\\
			\leq
			\left|v^{0}\right|\sqrt{\max\left(\left(\frac{4ad\theta}{\Delta x^{2}}+cd\theta -\lambda^{\alpha}d\theta-\left(1-b_{1}\right)\left(\theta+\left(1-\theta\right)e^{-\lambda\Delta t}\right)\right)^{2}, \left(cd\theta -\lambda^{\alpha}d\theta-\left(1-b_{1}\right)\left(\theta+\left(1-\theta\right)e^{-\lambda\Delta t}\right)\right)^{2}\right)+\left(\frac{bd\theta}{\Delta x}\right)^{2}}+\\
			\left|v^{0}\right| \sum_{j=1}^{k-1}\left(b_{j}-b_{j+1}\right)\left(\theta e^{-\lambda j \Delta t}+\left(1-\theta\right) e^{-\lambda \left(j+1\right) \Delta t}\right)+\left|v^{0}\right|\left|\theta e^{\lambda \Delta t}+1-\theta-d\lambda^{\alpha}-\sum_{j=0}^{k-1}\left(b_{j}-b_{j+1}\right)\left(\theta e^{-\lambda j \Delta t}+\left(1-\theta\right) e^{-\lambda \left(j+1\right) \Delta t}\right)\right|
			\\\leq
			\left\{ \sqrt{\max\left(\left(\frac{4ad\theta}{\Delta x^{2}}+cd\theta -\lambda^{\alpha}d\theta-\left(1-b_{1}\right)\left(\theta+\left(1-\theta\right)e^{-\lambda\Delta t}\right)\right)^{2}, \left(cd\theta -\lambda^{\alpha}d\theta-\left(1-b_{1}\right)\left(\theta+\left(1-\theta\right)e^{-\lambda\Delta t}\right)\right)^{2}\right)+\left(\frac{bd\theta}{\Delta x}\right)^{2}}+\right.\\
			\left.\sum_{j=1}^{N-1}\left(b_{j}-b_{j+1}\right)\left(\theta e^{-\lambda j \Delta t}+\left(1-\theta\right) e^{-\lambda \left(j+1\right) \Delta t}\right)+
			\left|\theta e^{\lambda \Delta t}+1-\theta-d\lambda^{\alpha}-\sum_{j=0}^{N-1}\left(b_{j}-b_{j+1}\right)\left(\theta e^{-\lambda j \Delta t}+\left(1-\theta\right) e^{-\lambda \left(j+1\right) \Delta t}\right)\right|\right\}\left|v^{0}\right|, \label{28}             \end{multline} 
		where the last inequality is true by (\ref{23}). By (\ref{26}), (\ref{27}) and (\ref{28}), equation (\ref{22}) implies \[\left|v^{0}\right|\geq\left|v^{k+1}\right|.\]
	\end{enumerate}
\end{Proof}\\
Note that for $\theta=0$ (\ref{22}) has a much simpler form. Moreover, only for that case, the stability can be provided independently on $\Delta x$.
\begin{Proposition}\label{oby2Prop}
	Let us assume
	\begin{equation}
		1-d\lambda ^{\alpha}-\displaystyle\sum_{j=0}^{N-1}\left(b_{j}-b_{j+1}\right)e^{-\lambda \Delta t (j+1)}\geq 0.\label{oby1}
	\end{equation}
	Then, for $\theta=0$ the scheme (\ref{17}) is stable.
\end{Proposition}
\begin{Proof}
	By (\ref{oby1}), the left side of (\ref{22}) for $\theta=0$ has the form \[cd-d\lambda^{\alpha}+1.\]
	On the other hand, the right side of (\ref{22}) is equal \[1-d\lambda^{\alpha}.\]
	Thus, (\ref{22}) holds. 	
\end{Proof}



\subsection{Convergence of the weighted discrete scheme}
We will show the convergence using the same techniques as in the case of stability. Let us define the error at the point $(x_{l},t_{k})$ by
$E_{l}^{k}=u_{l}^{k}-\hat{u}_{l}^{k}$, $l=0,1,\dots n$, $k=0,1,..,N$. Since $u$ and $\hat{u}$ have the same values in the initial condition, we conclude that $E_{i}^{0}=0$. Moreover, we assume that $x_{min}$ and $x_{max}$ are such that the error at the boundary is negligible, that is, $E_{0}^{k}=E_{n}^{k}=0$. Then, similarly to (\ref{19}) we get the following system:
\begin{equation}
	\begin{cases}
		-\left(\displaystyle\frac{ad}{\Delta x^{2}}+\displaystyle\frac{bd}{2\Delta x}\right)(\theta E^{0}_{i+1}+(1-\theta) E^{1}_{i+1})+
		\left(\displaystyle\frac{2ad}{\Delta x^{2}}+cd-d\lambda^{\alpha}\right)(\theta E^{0}_{i}+(1-\theta) E^{1}_{i})-\\ 
		\left(\displaystyle\frac{ad}{\Delta x^{2}}-\displaystyle\frac{bd}{2\Delta x}\right)(\theta E^{0}_{i-1}+(1-\theta) E^{1}_{i-1})= 
		(\theta e^{\lambda \Delta t}+(1-\theta)-d\lambda^{\alpha})E_{i}^{0}-(\theta e^{\lambda \Delta t}+1-\theta)E_{i}^{1}+
		\theta R_{i}^{0}+\left(1-\theta\right)R_{i}^{1},\\
		-\left(\displaystyle\frac{ad}{\Delta x^{2}}+\displaystyle\frac{bd}{2\Delta x}\right)(\theta E^{k}_{i+1}+
		(1-\theta) E^{k+1}_{i+1})+\left(\displaystyle\frac{2ad}{\Delta x^{2}}+cd-d\lambda^{\alpha}\right)(\theta E^{k}_{i}+(1-\theta) E^{k+1}_{i})-\\
		\left(\displaystyle\frac{ad}{\Delta x^{2}}-\displaystyle\frac{bd}{2\Delta x}\right)(\theta E^{k}_{i-1}+(1-\theta) E^{k+1}_{i-1})=
		\displaystyle\sum_{j=0}^{k-1}\left(b_{j}-b_{j+1}\right)(\theta e^{-j\Delta t\lambda}+(1-\theta) e^{-\left(j+1\right)\Delta t\lambda})(E_{i}^{k-j}-
		E_{i}^{0})+\\(\theta e^{\lambda \Delta t}+(1-\theta)-d\lambda^{\alpha})E_{i}^{0}-
		(\theta e^{\lambda \Delta t}+1-\theta)E^{k+1}_{i}+\theta R_{i}^{k}+\left(1-\theta\right)R_{i}^{k+1},\\
		E_{i}^{0}=0,\\
		E_{0}^{m}=E_{n}^{m}=0,
	\end{cases}\label{29}     
\end{equation}
where $k=1,\dots ,N-1$, $i=1,2,\dots ,n-1$, $m=0,\ldots,N$, and $\displaystyle\theta\in\left[0,1\right]$.
Similarly as in the case of stability, we will proceed with the von Neumann method. We introduce the following grid functions:\\
$E^{k}(x)=$
\begin{equation*}
	\begin{cases}
		\displaystyle E_{l}^{k}\text{, }x\in\left(x_{l}-\Delta x/2,x_{l}+\Delta x/2\right],\text{ }l=1,2\dots ,n-1,\\
		\displaystyle 0\text{, }x\in\left(x_{min},x_{min}+\Delta x/2\right]\cup\left[x_{max}-\Delta x/2,x_{max}\right].
	\end{cases}
\end{equation*}

$R^{k}(x)=$
\begin{equation*}
	\begin{cases}
		\displaystyle R_{l}^{k}\text{, }x\in\left(x_{l}-\Delta x/2,x_{l}+\Delta x/2\right],\text{ }l=1,2\dots ,n-1,\\
		\displaystyle 0\text{, }x\in\left(x_{min},x_{min}+\Delta x/2\right]\cup\left[x_{max}-\Delta x/2,x_{max}\right].
	\end{cases}
\end{equation*}

Because $E^{k}_{0}=E^{k}_{n}$, we make a periodic expansion of $E^{k}_{l}$ with the period $Y=x_{max}-x_{min}$. 
Then $E^{k}(x)$ has the following Fourier series extension:
\begin{equation*}
	\displaystyle E^{k}(x)=\sum_{j=-\infty}^{\infty}w_{j}^{k}e^{2j\pi x \iota/Y},
\end{equation*}
where $\displaystyle w_{j}^{k}=\frac{1}{Y}\int_{0}^{Y}E^{k}(x)e^{2j\pi x \iota/Y}dx$, $\iota=\sqrt{-1}$, $k=0,1,\dots N$. By analogy, because $R^{k}_{0}=R^{k}_{n}$, we make a periodic expansion for $R^{k}_{l}$ with the period $Y$. 
Then $R^{k}(x)$ has the following Fourier series extension:
\begin{equation*}
	\displaystyle R^{k}(x)=\sum_{j=-\infty}^{\infty}\rho_{j}^{k}e^{2j\pi x \iota/Y},
\end{equation*}
where $\displaystyle \rho_{j}^{k}=\frac{1}{Y}\int_{0}^{Y}R^{k}(x)e^{2j\pi x \iota/Y}dx$, $\iota=\sqrt{-1}$, $k=0,1,\dots N$. We define the norm $\norm{\cdot}_{\Delta x}$ as \[\displaystyle \norm{E^{k}}_{\Delta x}=\sqrt{\sum_{j=1}^{n-1}\Delta x\left|E_{j}^{k}\right|^{2}},\] \[\norm{R^{k}}_{\Delta x}=\sqrt{\sum_{j=1}^{n-1}\Delta x\left|R_{j}^{k}\right|^{2}},\] where 
\[\displaystyle E^{k}=\left(E_{1}^{k},E^{k}_{2},\dots ,E_{n-1}^{k}\right), \]
\[\displaystyle R^{k}=\left(R_{1}^{k},R^{k}_{2},\dots ,R_{n-1}^{k}\right). \]
Because $E_{0}^{k}=E_{n}^{k}=0,$ and $R_{0}^{k}=R_{n}^{k}=0$, there holds

\[\displaystyle \norm{E^{k}}_{\Delta x}^{2}=\int_{0}^{Y}\left|E^{k}\left(x\right)\right|^{2}dx=\norm{E^{k}\left(x\right)}^{2}_{L^{2}},\]
\[\displaystyle \norm{R^{k}}_{\Delta x}^{2}=\int_{0}^{Y}\left|R^{k}\left(x\right)\right|^{2}dx=\norm{R^{k}\left(x\right)}^{2}_{L^{2}}.\]
Using the Parseval identity we have:
\begin{equation}
	\begin{cases}
		\norm{E^{k}}_{\Delta x}^{2}=\displaystyle\sum_{j=1}^{n-1}\Delta x\left|E_{j}^{k}\right|^{2}=Y\sum_{j=-\infty}^{\infty}\left|w_{j}^{k}\right|^{2},\\
		\norm{R^{k}}_{\Delta x}^{2}=\displaystyle\sum_{j=1}^{n-1}\Delta x\left|R_{j}^{k}\right|^{2}=Y\sum_{j=-\infty}^{\infty}\left|\rho_{j}^{k}\right|^{2},
	\end{cases}\label{30}     
\end{equation}
where $k=0,1,\dots ,N.$ Based on the above analysis and the fact that $x_{l}=x_{min}+l\Delta x$ for $l=0,\ldots,n$, we 
suppose that the solution of (\ref{29}),
has the form:
\[E_{l}^{k}=w^{k}e^{\iota\beta \left(x_{min}+l\Delta x\right)},\]
\[R_{l}^{k}=\rho^{k}e^{\iota\beta \left(x_{min}+l\Delta x\right)},\]
where $\displaystyle\beta=\frac{2\pi j}{Y}.$
Substituting into (\ref{29}) we get:

\begin{equation}
	\begin{cases}
		\begin{gathered}
			\displaystyle \left(-\left(\frac{ad}{\Delta x^{2}}+\frac{bd}{2\Delta x}\right)
			e^{\iota\beta \Delta x}+\left(2\frac{ad}{\Delta x^{2}}+cd-d\lambda^{\alpha}\right)-\left(\frac{ad}{\Delta x^{2}}-\frac{bd}{2\Delta x}\right)
			e^{-\iota\beta \Delta x}\right)\left(\theta w^{0}+\left(1-\theta\right)w^{1}\right)=\\
			\displaystyle \left(\theta e^{\lambda \Delta t}+1-\theta-d\lambda^{\alpha}\right)w^{0}-\left(\theta e^{\lambda \Delta t}+1-\theta\right)w^{1}+\theta \rho^{0}+\left(1-\theta\right)\rho^{1},\\
			\displaystyle\left(-\left(\frac{ad}{\Delta x^{2}}
			+\frac{bd}{2\Delta x}\right)e^{\iota\beta \Delta x}+
			\left(2\frac{ad}{\Delta x^{2}}+cd-d\lambda^{\alpha}\right)
			-\left(\frac{ad}{\Delta x^{2}}-\frac{bd}{2\Delta x}\right)
			e^{-\iota\beta \Delta x}\right)
			\left(\theta w^{k}+\left(1-\theta\right)w^{k+1}\right)=\\
			\displaystyle\sum_{j=0}^{k-1}\left(b_{j}-b_{j+1}\right)\left( \theta e^{-j\Delta t\lambda}+\left(1-\theta\right)e^{-\left(j+1\right)\Delta t\lambda} \right)\left(w^{k-j}-w^{0}\right)+
			\left(\theta e^{\lambda \Delta t}+(1-\theta)-d\lambda^{\alpha}\right)w^{0}-\\
			(\theta e^{\lambda \Delta t}+1-\theta)w^{k+1}
			+\theta \rho^{k}+\left(1-\theta\right)\rho^{k+1},
		\end{gathered}
	\end{cases}\label{31}     
\end{equation}
where $k=1,\dots ,N-1$. By the identity \[ \displaystyle\sin^{2}\frac{z}{2}=-\frac{1}{4}\left(e^{\iota z}-2+e^{-\iota z}\right)\] 
and by $\rho^{0}=0$, $w^{0}=0$, we find that (\ref{31}) has the following form

\begin{equation}
	\begin{cases}
		\displaystyle \left(\zeta-d\lambda^{\alpha}\right)\left(1-\theta\right)w^{1}=-\left(\theta e^{\lambda \Delta t}+1-\theta\right)w^{1}
		+\left(1-\theta\right)\rho^{1},\\
		
		\displaystyle\left(\zeta-d\lambda^{\alpha}\right)\left(\left(1-\theta\right) w^{k+1}+\theta w^{k}\right)=
		\displaystyle\sum_{j=0}^{k-1}\left(b_{j}-b_{j+1}\right)\left( \theta e^{-j\Delta t\lambda}+\left(1-\theta\right)e^{-\left(j+1\right)\Delta t\lambda} \right)w^{k-j}-\\
		(\theta e^{\lambda \Delta t}+1-\theta)w^{k+1}
		+\theta \rho^{k}+\left(1-\theta\right)\rho^{k+1},
	\end{cases}\label{32}     
\end{equation}
where $k=1,\dots ,N-1$ and $\zeta$ were previously defined.
\begin{Lemma}\label{Lemmatemp}Let us assume that $\theta=0$ and
	\begin{equation}
		2d\left(c-\lambda^{\alpha}\right)+1\geq\sum_{j=0}^{N-1} (b_{j}-b_{j+1})e^{-\lambda \Delta t (j+1)}.\label{konieczartow1}
	\end{equation}
	Then:
	\begin{enumerate}
		\item if $\lambda^{\alpha}\neq c$, then $w^{k+1}$ follows 
		\[\displaystyle\left|w^{k+1}\right|\leq \frac{C_{2}\left|\rho^{1}\right|}{\left|c-\lambda^{\alpha}\right|d},\]
		where $k=0,1,\dots ,N-1$ and the constant $C_{2}$ is independent of $\Delta t$ and $\Delta x$.
		\item if $\lambda^{\alpha}=c$, then
		$w^{k+1}$ follows 
		\[\displaystyle\left|w^{k+1}\right|\leq \frac{C_{2}\left|\rho^{1}\right|}{b_{k}},\]
		where $k=0,1,\dots ,N-1$ and the constant $C_{2}$ is independent of $\Delta t$ and $\Delta x$.
	\end{enumerate}
\end{Lemma}

\begin{Proof}
	At the beginning, let us observe that assuming $\theta=0$ and (\ref{konieczartow1}), the left side of (\ref{22}) is equal to $1-d\lambda^{\alpha}+cd.$ Furthermore, by (\ref{konieczartow1}) we have $1-d\lambda^{\alpha}+cd>0$.
	Thus, by (\ref{26}) we have 	\begin{equation}
		\left|\zeta-d\lambda^{\alpha}+1\right| \geq 1-d\lambda^{\alpha}+cd>0.\label{wazne1}
	\end{equation}
	The convergence of the right series in the second line of (\ref{30}) implies that
	\begin{equation}\displaystyle\left|\rho^{k}\right|\equiv\left|\rho_{l}^{k}\right|\leq  C_{3}\left|\rho_{l}^{1}\right|\equiv C_{3}\left|\rho^{1}\right|,\label{zbych}
	\end{equation}
	for some positive constant $C_{3}$ and $k=2,\ldots,N$, $l=1,\ldots,n-1$
	Let us denote $C_{2}=\max(1,C_{3})$.
	\begin{enumerate}
		\item Let us observe that by (\ref{konieczartow1}) we have 
		\begin{equation}
			1+cd-d\lambda^{\alpha}\geq\left|cd-d\lambda^{\alpha}\right|+\sum_{j=0}^{N-1} (b_{j}-b_{j+1})e^{-\lambda \Delta t (j+1)}.\label{Agatkajestsuper!!}
		\end{equation}
		
		By the first equation of (\ref{32}), 
		we have \begin{equation*}
			\displaystyle\left|w^{1}\right|= \frac{\left|\rho^{1}\right|}{\left|\zeta-d\lambda^{\alpha}+1\right|}\leq
			\frac{\left|\rho^{1}\right|}{cd-d\lambda^{\alpha}+1}\leq \frac{\left|\rho^{1}\right|}{\left|cd-d\lambda^{\alpha}\right|}\leq
			\frac{C_{2}\left|\rho^{1}\right|}{d\left|c-\lambda^{\alpha}\right|}, \end{equation*} 
		where the first inequality holds by (\ref{wazne1}) and the second by (\ref{Agatkajestsuper!!}). Now let us suppose that    \begin{equation}
			\displaystyle \left|w^{m}\right|\leq 	\frac{C_{2}\left|\rho^{1}\right|}{d\left|c-\lambda^{\alpha}\right|},\label{35}      \end{equation} 
		where $m=1,2,\dots ,k$, $k<N$. By the second equation of (\ref{32}) and by (\ref{zbych}),
		we have \begin{multline}
			\displaystyle\left|\left(\zeta-d\lambda^{\alpha}+1\right)\right| \left|w^{k+1}\right|\leq
			\sum_{j=0}^{k-1}\left(b_{j}-b_{j+1}\right)e^{-\lambda \left(j+1\right)\Delta t }\left|w^{k-j}\right|+
			\left| \rho^{k+1}\right|\leq
			\sum_{j=0}^{k-1}\left(b_{j}-b_{j+1}\right)e^{-\lambda \left(j+1\right)\Delta t}\left|w^{k-j}\right|
			+C_{3}\left|\rho^{1}\right| \leq \\\ 
			\displaystyle\left(\sum_{j=0}^{N-1}\left(b_{j}-b_{j+1}\right)e^{-\lambda \left(j+1\right)\Delta t}+\left|cd-d\lambda^{\alpha}\right|\right)\frac{C_{2}\left|\rho^{1}\right|}{d\left|c-\lambda^{\alpha}\right|}\leq \left(1-d\lambda^{\alpha}+cd\right)\frac{C_{2}\left|\rho^{1}\right|}{d\left|c-\lambda^{\alpha}\right|},\label{36}     
		\end{multline}
		where the third inequality is true by (\ref{35}) and the last by (\ref{Agatkajestsuper!!}). Dividing (\ref{36}) by\\ $\left|\zeta-d\lambda^{\alpha}+1\right|$, by (\ref{wazne1}) we get
		\[ \displaystyle\left|w^{k+1}\right|\leq	\frac{C_{2}\left|\rho^{1}\right|}{d\left|c-\lambda^{\alpha}\right|}.\]
		
		By mathematical induction, the proof of the first part is completed.

		\item At the beginning, let us observe that by Proposition \ref{3.1} \begin{dmath}
			1\geq
			\sum_{j=0}^{N-1}\left(b_{j}-b_{j+1}\right) e^{-\lambda \left(j+1\right) \Delta t}+b_{N}.\label{22b}  
		\end{dmath} holds.
		By the first equation of (\ref{32}), 
		we have \begin{equation*}
			\displaystyle\left|w^{1}\right|= \frac{\left|\rho^{1}\right|}{\left|\zeta-d\lambda^{\alpha}+1\right|}\leq
			\frac{\left|\rho^{1}\right|}{1}\leq
			\frac{C_{2}\left|\rho^{1}\right|}{b_{0}}, \end{equation*} 
		where the first inequality holds by (\ref{wazne1}) and by $c=\lambda^{\alpha}$. Now let us suppose that    \begin{equation}
			\displaystyle \left|w^{m}\right|\leq 	\frac{C_{2}\left|\rho^{1}\right|}{b_{k-1}},\label{35a}      \end{equation} 
		where $m=1,2,\dots ,k$, $k<N$. By the second equation of (\ref{32}) and by (\ref{zbych}),
		we have \begin{multline}
			\left|w^{k+1}\right|=\left|\sum_{j=0}^{k-1}\left(b_{j}-b_{j+1}\right)\left(e^{-\lambda \left(j+1\right)\Delta t }\right)w^{k-j}+
			\rho^{k+1}\right|\leq 
			\sum_{j=0}^{k-1}\left(b_{j}-b_{j+1}\right)e^{-\lambda \left(j+1\right)\Delta t}\left|w^{k-j}\right|
			+C_{3} \left|\rho^{1}\right| \leq \\\ 
			\left(\sum_{j=0}^{k-1}\left(b_{j}-b_{j+1}\right)e^{-\lambda \left(j+1\right)\Delta t}+b_{k}\right)
			\frac{C_{2}\left|\rho^{1}\right|}{b_{k}}\leq \frac{C_{2}\left|\rho^{1}\right|}{b_{k}},\label{36a}     
		\end{multline}
		where the second inequality is true by (\ref{35a}) and the last by (\ref{22b}).
		By mathematical induction, the proof of the second part is completed.
	\end{enumerate}
\end{Proof}\\

\begin{Theorem}\label{kdom}
	Let us assume the conditions of Lemma \ref{Lemmatemp}, then the discrete scheme (\ref{17}) is convergent and it follows that
	\[\displaystyle\norm{u_{l}^{k}-\hat{u}_{l}^{k}}_{\Delta x}\leq  C_{4}  \left(\Delta t^{2-\alpha}+\Delta x^{2}\right),\]
	for $k=1,2,\dots ,N,$ where $C_{4}$ is a positive constant independent of $\Delta t$ and $\Delta x$.
\end{Theorem}

\begin{Proof}
	By Theorem \ref{rukuio} there exists a positive constant $C_{1}$, such that
	for $\Delta t$ and $\Delta x$ small enought there holds \[\left|R_{i}^{k}\right|\leq C_{1}\Delta t^{\alpha}\left(\Delta t^{2-\alpha}+\Delta x^{2}\right),\]		where $k=1,\ldots,N$, $j=1,\ldots,n-1$. Then by (\ref{30}) we get
	\begin{equation}
		\displaystyle\norm{R^{k}}_{\Delta x}\leq C_{1}\sqrt{Y}\Delta t^{\alpha}\left(\Delta t^{2-\alpha}+\Delta x^{2}\right).\label{33}     
	\end{equation}
	
	By Lemma \ref{Lemmatemp}, for $c\neq\lambda^{\alpha}$ and $k=1,\dots,N$ we have
	\begin{equation*}\displaystyle\left|w^{k}\right|\leq\frac{C_{2}\left|\rho^{1}\right|}{\left|c-\lambda^{\alpha}\right|d}.
	\end{equation*}
	
	Similarly, by (\ref{30}) and (\ref{33}), we have the following:
	\[\displaystyle\norm{E^{k}}_{\Delta x}\leq\frac{C_{2}}{d\left|c-\lambda^{\alpha}\right|}\norm{R^{1}}_{\Delta x}\leq \frac{C_{2}}{\Gamma(2-\alpha)\left|c-\lambda^{\alpha}\right|}C_{1}\sqrt{Y}\left(\Delta t^{2-\alpha}+\Delta x^{2}\right).\]
	Then $\displaystyle C_{4}=\frac{C_{2}}{\left|c-\lambda^{\alpha}\right|\Gamma(2-\alpha)}C_{1}\sqrt{Y}$. If $c=\lambda^{\alpha}$, then by Lemma \ref{Lemmatemp}, for $k=1,\ldots,N$ we have 
	\[\displaystyle \left|w^{k}\right|\leq\frac{C_{2}\left|\rho^{1}\right|}{b_{k-1}}\leq
	\frac{C_{2}\Delta t^{-\alpha} T^{\alpha}\left|\rho^{1}\right|}{b_{N-1}N^{\alpha}}\leq \frac{C_{2}\Delta t^{-\alpha} T^{\alpha}\left|\rho^{1}\right|}{1-\alpha},\]
	where the second inequality holds by Proposition \ref{3.1} and the last by (\ref{1-alf}). Similarly
	\[\norm{E^{k}}_{\Delta x}\leq \frac{C_{2} \Delta t^{-\alpha}T^{\alpha} \norm{R^{1}}_{\Delta x}}{1-\alpha}\leq \frac{C_{2} T^{\alpha} }{1-\alpha}C_{1}\sqrt{Y}\left(\Delta t^{2-\alpha}+\Delta x^{2}\right).\]
	Then $\displaystyle C_{4}=\frac{C_{2}T^{\alpha}}{1-\alpha}C_{1}\sqrt{Y}$. \end{Proof}

\begin{figure}[ht]
	\centering
	\includegraphics[scale=0.34]{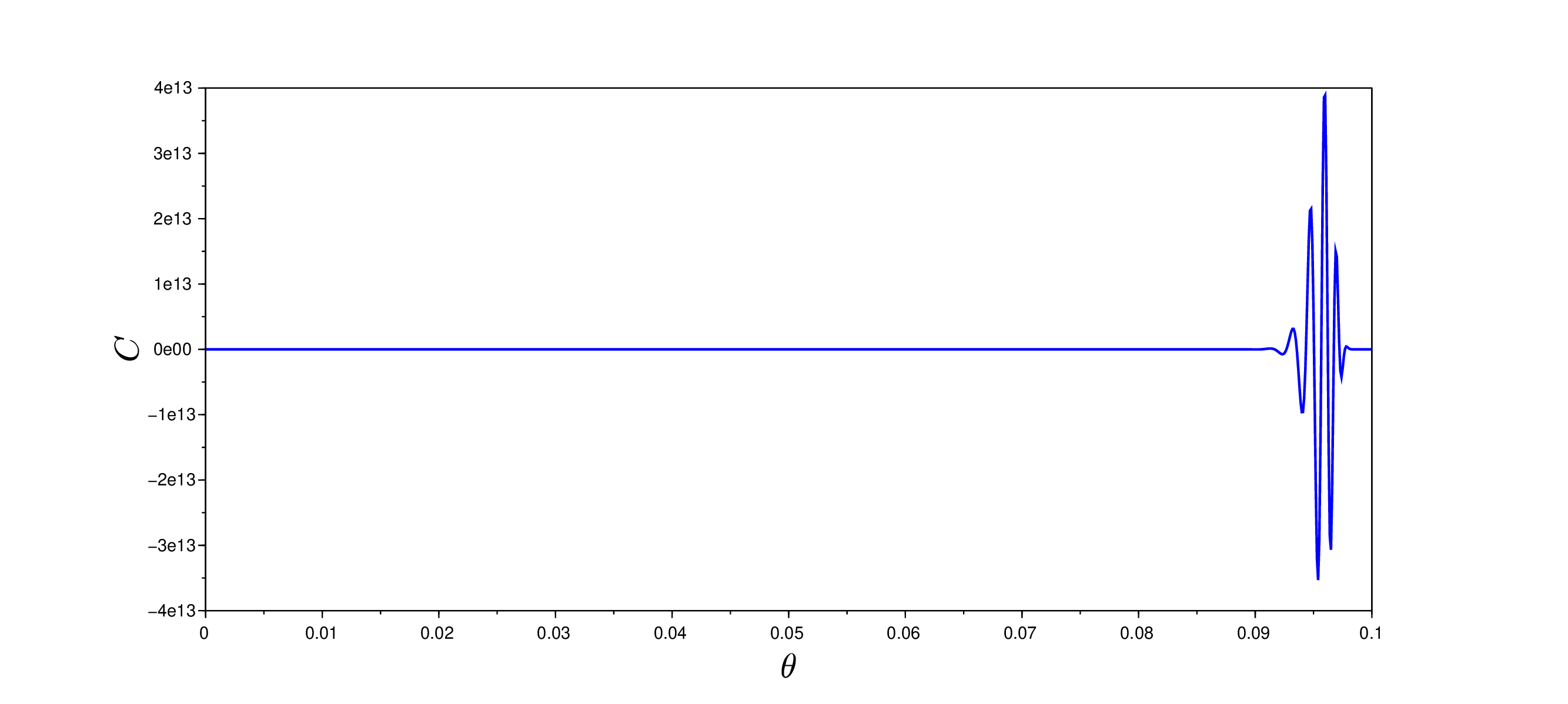}
	\caption{The European call price in dependence of $\theta$. The result of lack of stability is observed. The parameters: $T=\sigma=\lambda=0.3$, $S_{0}=1$, $r=0.08$, $x_{min}=-20$, $x_{max}=10$, $n=5000$, $N=100$,
		$\alpha=0.8$, $K=2$,
		.}\label{temptheta}
\end{figure}

We recall the observation of \cite{heston2000rate} that the rate of convergence depends on the smoothness of the option payoff
function. Since the payoff function is not continuously differentiable (in the case of the European option, the critical point is at the strike), the rate of convergence can be lower than indicates the previous theorem. This problem can resolve the smoothing of the payoff function, for example, by transforming the original payoff function $f(x)$ into $f^{\ast}(x)=\frac{1}{2\Delta x}\int_{-\Delta x}^{\Delta x} f(x-y)dy$ \cite{heston2000rate}. The same remark is true for both the finite difference (FD) and the CRR method (the method considered in the numerical examples).
A further analysis of the impact of smoothing the payoff function extends the scope of this paper. 
Let us observe that because of the lack of unconditional stability/convergence (i.e., the stability/convergence which is independent of $\Delta t$ and $\Delta x$) we can not decide which value of the parameter $\theta$ is the most optimal. In \cite{ja2} it was shown that for the subdiffusive B-S model (case $\lambda=0$) the most optimal value is $\displaystyle\check{\theta_{\alpha}}=\frac{2-2^{1-\alpha}}{3-2^{1-\alpha}}$. Then the lowest boundary for an error is achieved with the conservation of unconditional stability/convergence. We cannot repeat this approach since (\ref{17}) is not unconditionally stable/convergent. However, for $\theta=0$ and $\lambda^{\alpha}\leq c$, (\ref{konieczartow1}) holds for all $\Delta t$ and $\Delta x$. Thus, by Theorem \ref{kdom} in that particular case the numerical scheme (\ref{17}) is unconditionally convergent. 
Note that all considerations - especially (\ref{22}) - are much easier if $\theta=0$. Only in this scenario, the stability/convergence of (\ref{17}) can be obtained independently on $\Delta x$. Moreover, we have shown the conditional convergence only for $\theta=0$. In our opinion, the implicit scheme (i.e., the scheme (\ref{17}) for $\theta=0$) has the most practical impact. In Figure \ref{temptheta} the relation between the price of European call option and the parameter $\theta$ is presented. For the considered parameters and $\theta=0$, (\ref{konieczartow1}) is satisfied. However, with increasing $\theta$, the numerical scheme (\ref{17}) loses stability. Note that in the case of the subdiffusive B-S model, the unconditional stability/convergence is provided for $\theta\in[0,\displaystyle\check{\theta_{\alpha}}]$, and for $\theta\in(\displaystyle\check{\theta_{\alpha}},1]$ there is a conditional stability/convergence \cite{ja2}. In the classical B-S the same property is conserved because $\displaystyle\check{\theta_{1}}=1/2$ \cite{ja}. For the case considered in Figure \ref{temptheta} $\displaystyle\check{\theta_{\alpha}}=0.46$. We can conclude that each generalization of the B-S model is paid by a more narrow interval for $\theta$ where the unconditional stability/convergence holds. We can also observe that with the generalization of B-S model, the conditions providing convergence are getting more complex.
Note that for given parameters $T$, $\alpha$, $\lambda$, $r$ finding $N$ such that (\ref{konieczartow1}) is satisfied can be difficult or even impossible, especially for the "big" $\lambda$.

Let us observe that the tsB-S model (similarly to a standard B-S model) uses the expiration time in a non-dimensional form - in other words, $T=1$ can be $1$ day or $1$ year. It is important to highlight that an interest rate $r$, a dividend rate $\delta$, volatility $\sigma$ and tempered subdiffusion parameters $\alpha$, $\lambda$ have to be related with the considered unit of time. This is the essence of the following fact:
	\begin{Proposition}
		Let us denote $C_{tsB-S}(Z_{0},K,T,r,\delta,\sigma,\alpha,\lambda)$ as a fair price of the European call option in the tsB-S model. Then, for $\beta>0$ we have:
		\begin{equation}
			C_{tsB-S}(Z_{0},K,T,r,\delta,\sigma,\alpha,\lambda)=C_{tsB-S}(Z_{0},K,\beta T,\frac{r}{\beta^{\alpha}},\frac{\delta}{\beta^{\alpha}},\frac{\sigma}{\sqrt{\beta^{\alpha}}},\alpha,\frac{\lambda}{\beta}).\label{skalow_propert}
		\end{equation}	
\end{Proposition}
\begin{Proof}
	Let us consider a change of scale of time:
		\begin{equation}
			t^{*}=\beta t,\label{eurtemptranstemp}
		\end{equation} for $\beta>0$ and $t\in[0,T]$.
		An interest rate after transformation (\ref{eurtemptranstemp}) is related to the old one by
		\[e^{r^{*}T^{*}}=e^{rT},\]
		Hence, \begin{equation}
			r^{*}=\frac{r}{\beta}.\label{eurtemptranstemp93}
		\end{equation}
		After the time transformation (\ref{eurtemptranstemp}), $\sigma$ will change as follows:
		\begin{equation}
			\sigma^{*}=\frac{\sigma}{\sqrt{\beta}}.\label{eurtemptranstemp94}
		\end{equation}
		Therefore, based on (\ref{eurtemptranstemp})-(\ref{eurtemptranstemp93}) and (\ref{eurtemptranstemp94}) we get
		\begin{equation}
			C(Z_{0},K,T,r,\delta,\sigma)=C(Z_{0},K,T\beta,\frac{r}{\beta},\frac{\delta}{\beta},\frac{\sigma}{\sqrt{\beta}}),\label{rtlzuao4}
		\end{equation}  where $\beta>0$ and $C(Z_{0},K,T,r,\delta,\sigma)$ denotes a fair price of European call option in the standard B-S model. It is worth to observe, that (\ref{rtlzuao4}) can be easily verified. We know, that the fair price of European call option in dependence on $(Z_{0},K,T\beta,\frac{r}{\beta},\frac{\delta}{\beta},\frac{\sigma}{\sqrt{\beta}})$ (for $\beta>0$) is given by B-S formula \cite{McDonald}:
		\begin{equation}
			\displaystyle C(Z_{0},K,T\beta,\frac{r}{\beta},\frac{\delta}{\beta},\frac{\sigma}{\sqrt{\beta}})=Z_{0}e^{-\frac{\delta}{\beta} T\beta}\Phi(d1)-Ke^{-\frac{r}{\beta}T\beta}\Phi(d_{2}),\label{idide}
		\end{equation}
		where $	d_{1}=\displaystyle\frac{\log\frac{Z_{0}}{K}+(\frac{r}{\beta}-\frac{\delta}{\beta}+\frac{1}{2}(\frac{\sigma}{\sqrt{\beta}})^2)T\beta}{\frac{\sigma}{\sqrt{\beta}}\sqrt{T\beta}}$ and $
		\displaystyle d_{2}=d_{1}-\frac{\sigma}{\sqrt{\beta}}\sqrt{T\beta}$.
		Since $\beta$ in (\ref{idide})  can be reduced, so by (\ref{idide}) we get (\ref{rtlzuao4}).
By \cite{Meser} we have
		\begin{equation*}\vartheta_{\alpha,\lambda}(x,\beta t)=\frac{1}{\beta ^{\alpha}}\vartheta_{\alpha,\beta\lambda}\left(\frac{x}{\beta^{\alpha}}, t\right),\label{mizii}
		\end{equation*}
		which can be summarized as
		\begin{equation}S_{\alpha,\lambda}(\beta t)\stackrel{d}{=}\beta ^{\alpha}S_{\alpha,\beta\lambda}(t).\label{mizii2}
		\end{equation}
		Let us observe that for $\beta>0$ we have:
		\begin{multline}
			C_{tsB-S}(Z_{0},K,T\beta,r,\delta,\sigma,\alpha,\lambda)=\E C(Z_{0},K,S_{\alpha,\lambda}(T\beta),r,\delta,\sigma,\alpha,\lambda)=\\\E C(Z_{0},K,\beta^{\alpha}S_{\alpha,\beta\lambda}(T),r,\delta,\sigma)=\E C(Z_{0},K,S_{\alpha,\beta\lambda}(T),r\beta^{\alpha},\delta\beta^{\alpha},\sigma\sqrt{\beta^{\alpha}})=\\C_{tsB-S}(Z_{0},K,T,r\beta^{\alpha},\delta\beta^{\alpha},\sigma\sqrt{\beta^{\alpha}},\alpha,\beta\lambda),\label{wrocwie7}
		\end{multline}
		where the first and last equalities are satisfied by \cite{Gajda}, the second equality by (\ref{mizii2}), and the third by (\ref{rtlzuao4}). 
		Based on (\ref{wrocwie7}) the proof is completed.
\end{Proof}

\subsection{Numerical examples}
\begin{Example}
	The tsB-S is the generalization of the subdiffusive B-S model \cite{ja3, ja2, MM} therefore, we can expect that for $\lambda$ close to $0$ the prices in both models will be almost the same. It is important to note that for small $\alpha$ and $\lambda$ close to $0$, the term $\lambda^{\alpha}$ can be close to $1$. Then, even for $\lambda$ close to $0$, (\ref{12}) is not reducing to the subdiffusive PDE \cite{ja2} (to be so, $\lambda^{\alpha}$ has to diminish). Let us consider the following parameters: $\lambda=10^{-10}$, $T=1$, $K=2$, $Z_{0}=1$, $\theta=0$, $x_{min}=-10$,	$x_{max}=10$, $r=\sigma=0.5$, $n=N=900$. In Table \ref{table:1p1}, the prices of European call options in the tsB-S and subdiffusive models for different $\alpha$ are presented. The last column of the table is the relative difference, i.e., $|C_{tsB-S}-C_{subB-S}|/C_{subB-S}$. We see that for a small value of $\alpha$, even $\lambda$ close to $0$ cannot guarantee that both models will return almost the same result. However, with $\lambda\rightarrow0$, tsB-S will be reduced to the subdiffusive model. The essential question is how small $\lambda$ should be to allow the tsB-S model for a given $\alpha$ to be close to the subdiffusive B-S model. Consider the threshold to be equal to $1\%$. If such an effect is caused by the term $\lambda^{\alpha}$, we have to provide at least $\lambda^{\alpha}=1\%$. In Table \ref{table:1p2} the prices of European call options are presented in the tsB-S and subdiffusive B-S models for different parameters $\alpha$ and $\lambda=0.01^{1/\alpha}$. We see that for the assumed threshold $1\%$, the relative difference $|C_{tsB-S}-C_{subB-S}|/C_{subB-S}$ is between $1.33\%$ and $1.35\%$. Therefore, we conclude that the effect discussed is caused by the term $\lambda^{\alpha}$. Note that for the considered parameters and all $\alpha$ from Table \ref{table:1p1}, Table \ref{table:1p2} and Table \ref{table:2p2}, the convergence condition (\ref{konieczartow1}) is satisfied.
	\begin{table}[!ht]\footnotesize
		\begin{center}
			\begin{tabular}{|c c c c|} 
				\hline
				$\alpha$&   $C_{tsB-S}\times10^{2}$  &  $C_{subB-S}\times10^{2}$&relative difference\\ \hline
				$\alpha=10^{-11}$& $97.86$ &  $15.11$ &  $547.66\%$\\		
				$\alpha=0.01$& $58.86$ &  $15.22$ &  $286.69\%$\\		
				$\alpha=0.1$& $18.51$ &  $16.15$ &  $14.67\%$\\		
				$\alpha=0.2$& $17.2$ &  $16.97$ &  $1.35\%$\\		
				$\alpha=0.3$& $17.59$ &  $17.57$  & $0.13\%$\\		
				$\alpha=0.4$& $17.91$ &  $17.91$ &  $0.01\%$\\		
				$\alpha=0.5$& $17.97$ &  $17.97$ &  $1.21\times 10^{-5}$\\	
				$\alpha=0.6$& $17.7$ &  $17.7$ &   $1.1\times 10^{-6}$\\		
				$\alpha=0.7$& $17.06$ &  $17.06$ &   	$10^{-7}$\\
				$\alpha=0.8$& $16.03$ &  $16.03$ &   $9.3\times 10^{-9}$\\		
				$\alpha=0.9$& $14.58$ &  $14.58$ &   $7.58\times 10^{-10}$\\		
				$\alpha=0.999$& $12.77$ &  $12.77$ &   $1.82\times 10^{-12}$\\			\hline
			\end{tabular}
			\caption{\label{table:1p1}The prices of the European call option in tsB-S and subdiffusive B-S for $\lambda=10^{-10}$ and different $\alpha$. For small $\alpha$ the prices provided by both models are not close to each other.}
		\end{center}
	\end{table}

	\begin{table}[!ht]\footnotesize
		\begin{center}
			\begin{tabular}{|c c c c c|} 
				\hline
				$\alpha$&  $\lambda=0.01^{1/\alpha}$ & $C_{tsB-S}\times10^{2}$  &  $C_{subB-S}\times10^{2}$&relative difference\\ \hline
				$0.01$ & $ 10^{-200} $ &  $15.42$ &   $15.22$ &   $1.33\%$ \\
				$0.02$ &   $10^{-100} $ & $ 15.54$ &   $15.33$ &   $1.33\%$ \\
				$0.03$ &   $2.15\times10^{-67}$ &  $15.65$ &   $15.44$ &   $1.33\%$ \\
				$	0.04$ &   $10^{-50} $ &  $15.76$ &   $15.55$ &   $1.34\%$ \\
				$0.05$ &  $ 10^{-40}$ &  $ 15.86$ &  $ 15.65$ &   $1.34\%$ \\
				$0.06$ &  $ 4.64\times 10^{-34}$ &  $ 15.97$ &   $15.76$ &   $1.34\%$ \\
				$0.07$ &  $ 2.68\times10^{-29} $ &  $16.07$ &  $ 15.86$ &    $1.34\%$ \\
				$0.08$ &  $ 10^{-25}$ &   $16.17$ &  $ 15.96$ &   $1.34\%$ \\
				$0.09$ &  $ 6\times10^{-23} $ &  $16.27$ &  $ 16.05$ &   $1.35\%$ \\
				$0.1 $ &   $10^{-20} $ &  $16.36 $ & $ 16.15$ &   $1.35\%$
				\\			\hline
			\end{tabular}
			\caption{\label{table:1p2}The prices of the European call option in tsB-S and subdiffusive B-S models for different $\alpha$ and $\lambda=0.01^{1/\alpha}$. For each $\alpha$ the prices provided by both models are close to each other.}
		\end{center}
	\end{table}

	Finally, we expect that for $\alpha$ close to $1$, the prices of the option in tsB-S should be close to the prices provided by the subdiffusive B-S and B-S model. Let us take $\alpha=0.999$. In Table \ref{table:2p2} we investigate how close tsB-S and B-S are, depending on $\lambda$. The last row of the table is the relative difference, i.e., $|C_{tsB-S}-C_{B-S}|/C_{B-S}$. The related price of the European call option in the subdiffusive B-S and B-S models is $0.1277$ and $0.1276$, respectively. The models return results that are very close to each other. Note that as $\lambda$ increases, the relative difference also increases. 
	
	\begin{table}[!ht]\footnotesize
		\begin{center}
			\begin{tabular}{|c| c c c c c c |} 
				\hline
				$\lambda$&   $10^{-11}$  &  $0.01$&$0.1$&$0.4$&$0.7$&$0.9$\\ 
				$C_{tsB-S}\times10^{2}$ & $12.77$& $12.78$& $12.78$& $12.8$& $12.81$& $12.82$\\
				relative difference&$0.1\%$&$0.11\%$&$0.16\%$&$0.29\%$&$0.41\%$&$0.49\%$\\
				\hline
			\end{tabular}
			\caption{\label{table:2p2}The prices of European call options in the tsB-S model with a relative difference from the corresponding prices of the option in the B-S model for $\alpha=0.999$ and different $\lambda$. In each case, the prices provided by both models are close to each other.}
		\end{center}
	\end{table}

	In Figure \ref{Figure3} we present how the price of the European call option $C_{tsB-S}$ behaves for different values of $\alpha$ and $\lambda$. To correctly present the shape of $C_{tsB-S}(\alpha,\lambda)$, we have to use a more dense mesh for $\lambda$ close to $0$. 
	It is also worth mentioning that $C_{tsB-S}(\alpha,\lambda)$ is not an injective function. The parameters are $\sigma=r=T=Z_{0}=1$, $K=2$, $n=1000$, $N=400$, $x_{min}=-10$, $x_{max}=10$, $\theta=0$.

	\begin{figure}[ht]
		\centering
		\includegraphics[scale=0.44]{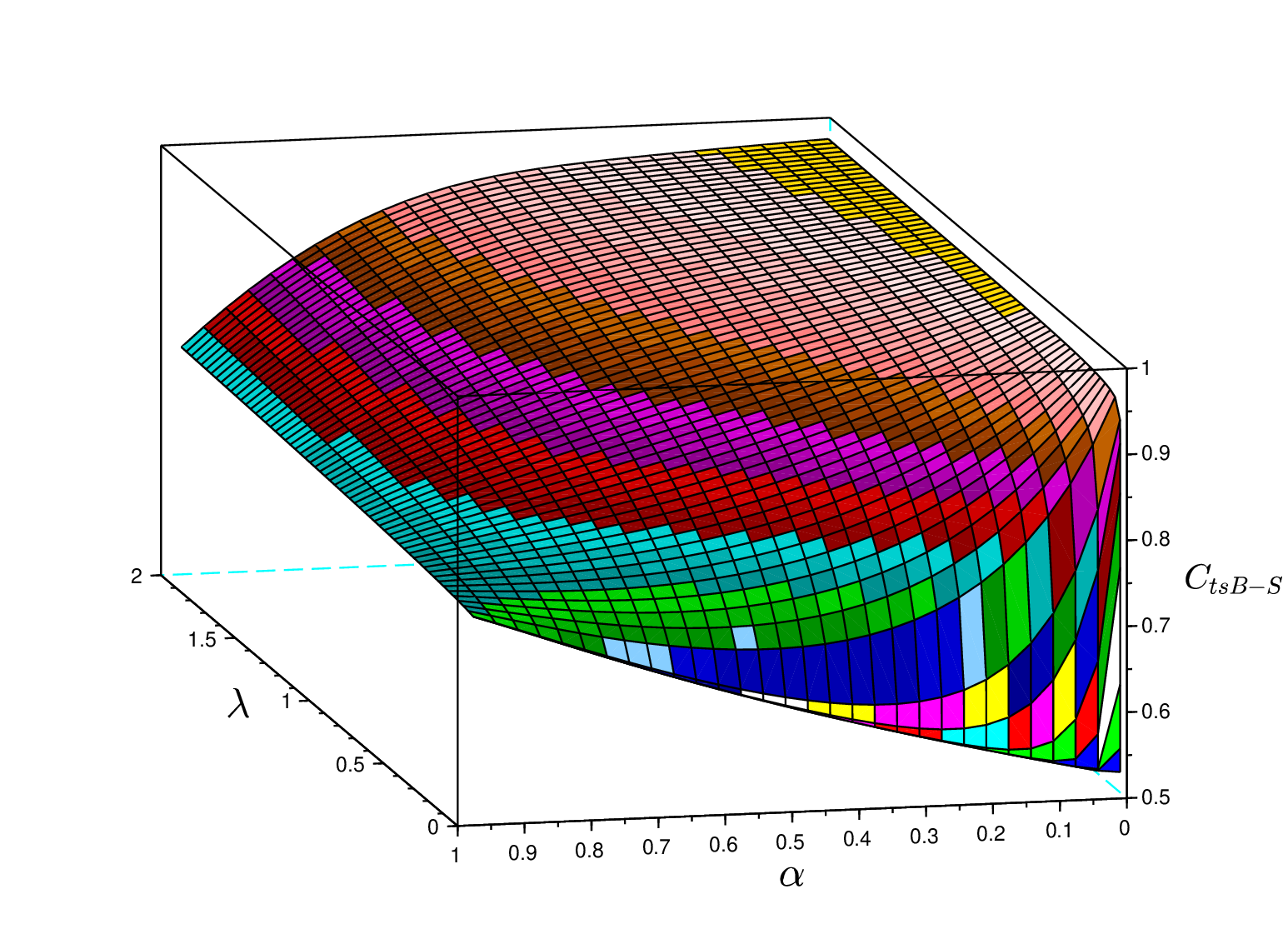}
		\caption{The dependence of the European call price on $\lambda$ and $\alpha$ for $T=1$. Note that for $\lambda=0$, the argument $T=1$ is the inflection point of the function $C_{tsB-S}(\alpha,T)$ \cite{ja2}. } \label{Figure3}  
	\end{figure}		
\end{Example}

\newpage
\begin{Example}\label{eweer}
	Let us take the parameters $Z_{0}=1$, $\sigma=1$, $\alpha=0,8$, $\theta=0$, $x_{max}=10$, $x_{min}=-20$, $n=70$, $N=360$.
	In Figure \ref{Figure2} we present the dependence of the fair price of the European call option on $K$ and $T$. We see how the well-known interactions from the classical and subdiffusive cases are conserved in tsB-S. For Figure \ref{Figure2} we take $r=\lambda=1$. Let us observe that for parameters of Figures \ref{Figure3}-\ref{Figure2}, (\ref{konieczartow1}) holds.
	\begin{figure}[ht!]
		\centering
		\includegraphics[scale=0.425]{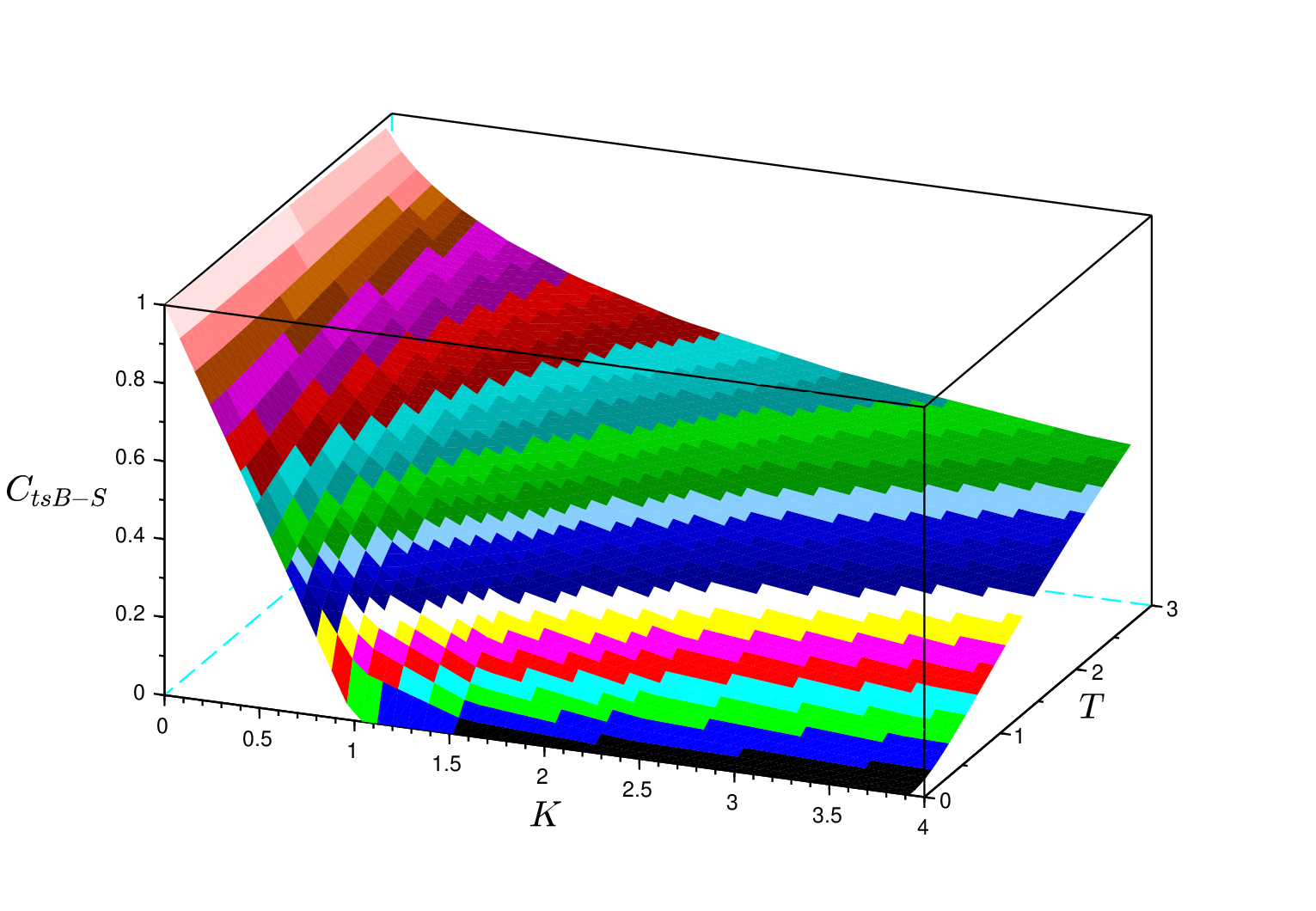}
		\caption{The dependence on $K$ and $T$ of the European call price in tsB-S model. The behavior of the option price in the dependence of financial parameters is analogous to the classical B-S model.} \label{Figure2}  
	\end{figure} 
	In Figures \ref{Figure4} (a) and \ref{Figure4} (b) we compare the FD with the MC (Monte Carlo) method (for $M_{1}=400$ repetitions and $k_{1}=50$ intermediate points) explained in \cite{Gajda} and the CRR method explained in \cite{jaCRR} (for $M_{2}=400$ repetitions and $k_{2}=40$ intermediate points).
	The other parameters are
	$N=360$, $K=2$. For such parameters, all
	methods compute the price of the option at approximately the same time $1,1s$. We conclude that MC and CRR are
	less accurate than FD output. Increasing the value of numerical parameters $k$ and $M$ of both methods
	follows their
	output is approaching the real fair price but also increases the time of computation. Although for most of the values $\lambda$
	in Figure \ref{Figure4} (a) the convergence condition (\ref{konieczartow1}) is not satisfied, the output of the FD method is in the regime indicated by the MC and CRR methods. However, if (\ref{konieczartow1}) is not provided, the output of the FD method may be different from the real value of the option. Particular caution should be taken for the "big" $\lambda$. In Figure \ref{Figure4} (b) the effect of lack of convergence is observed. The best practice is to use the method if its convergence conditions are provided.\newpage 
	
	\begin{figure}[ht!]
		\begin{center}
			\subfigure[]{
				\includegraphics[scale=0.32]{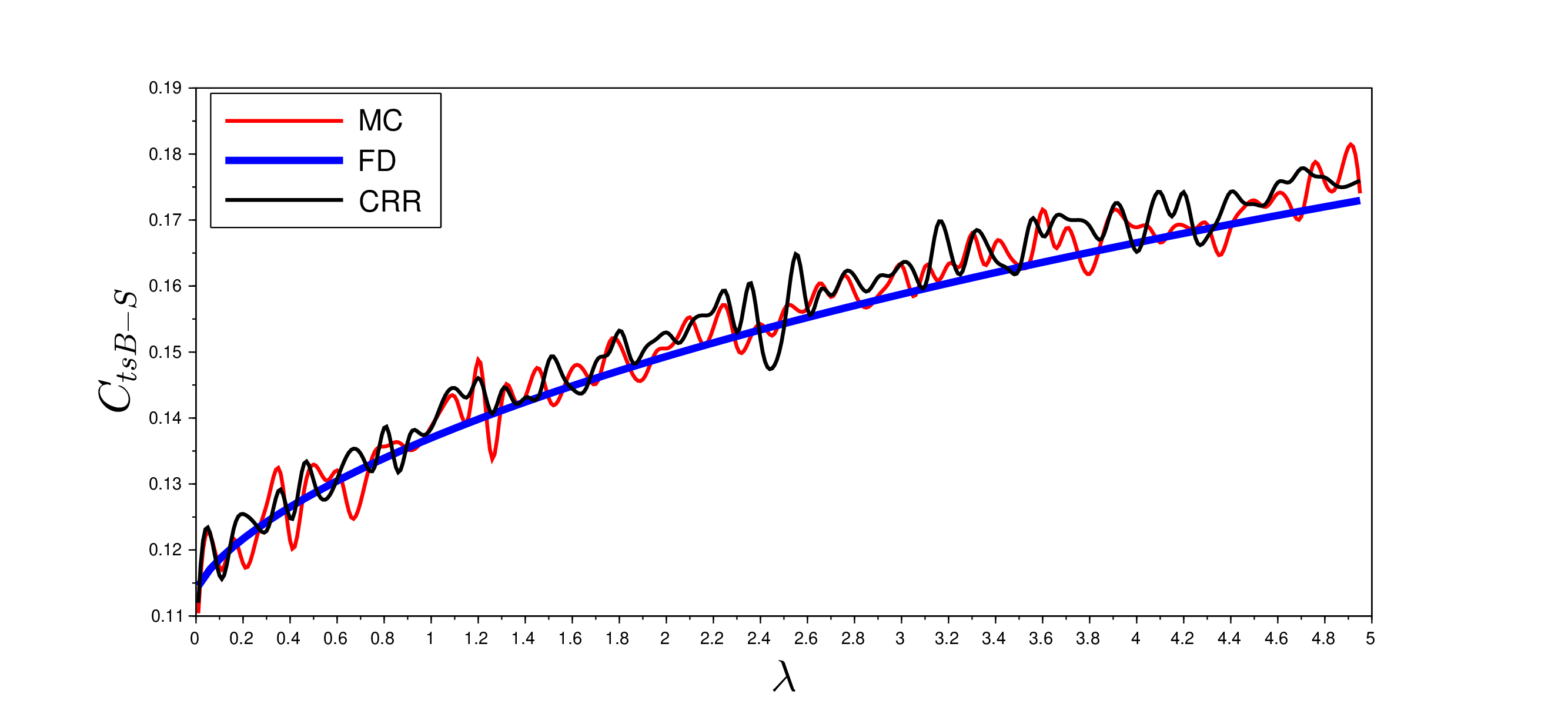}
			}
			\subfigure[]{
				\includegraphics[scale=0.32]{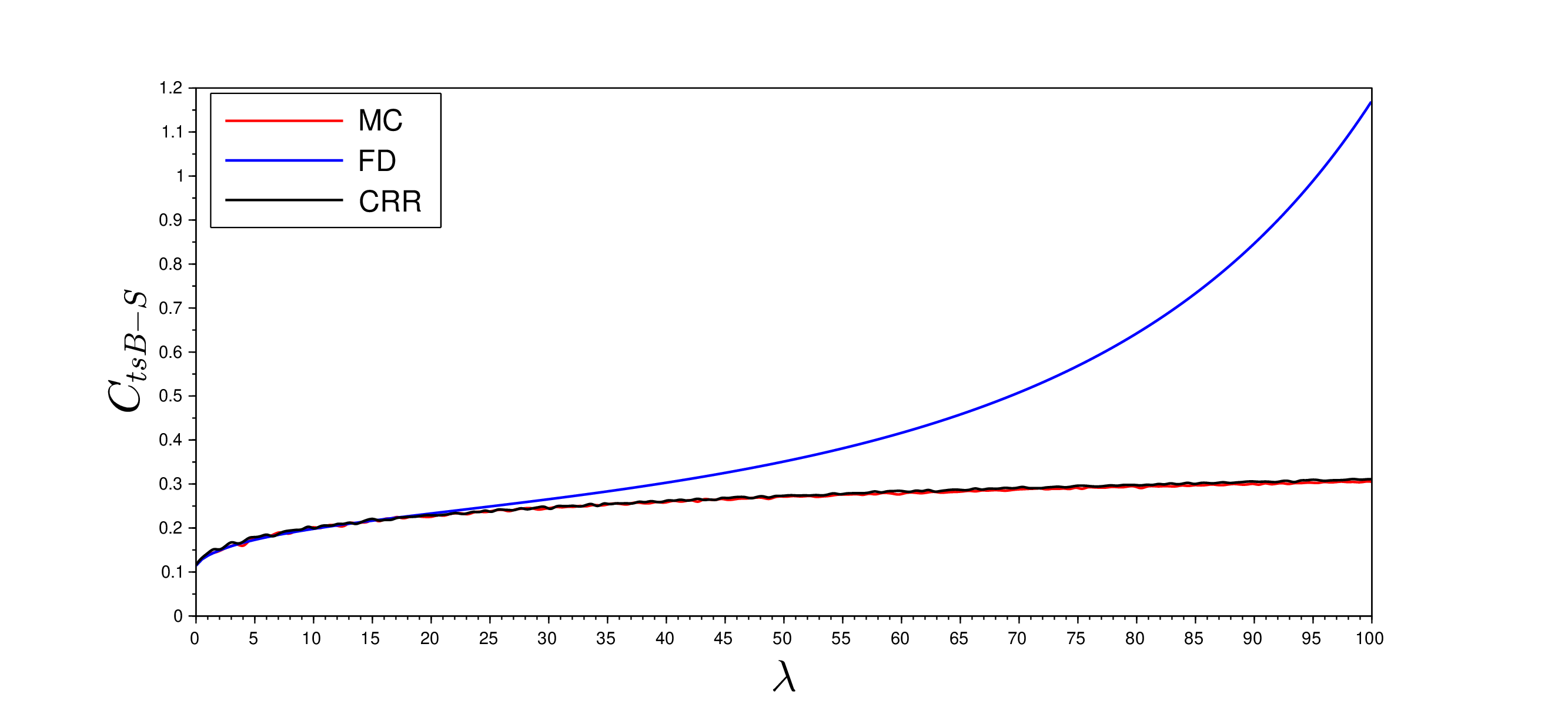}
			}
		\end{center}
		\caption{The price of European call option in dependence of $\lambda$ for $T=0.5$. }   \label{Figure4}  
	\end{figure}  
\end{Example} 
	\newpage\begin{Example}
		Let us consider the parameters:
		$T=K=\sigma=r=\lambda=1$, $\delta=0$, $Z_{0}=2$, $n_{FD}=200$, $N_{FD}=200$, $\theta=0$, $x_{max}=10$, $x_{min}=-20$. The equality (\ref{skalow_propert}) can be used to optimize numerical methods. Tables \ref{FDblad:1p3} - \ref{CRRblad:1p3} present relative errors (multiplied by $10^{3}$) calculated by $\alpha$ and scale parameter $\beta$, for the FD, MC and CRR methods (the latter two were introduced in Example \ref{eweer}). MC and CRR methods were used for $M_{MC}=M_{CRR}=1000$ repetitions and
		$k_{CRR}=k_{MC}=50$ intermediate points. Furthermore, for the CRR method, we take $n_{CRR}=50$, where $n_{CRR}$ is the number of nodes in the classical CRR method. The exact value is approximated by the result of the FD method for $N_{FD}=800$, $n_{FD}=300$, and $\beta=1$. Tables \ref{FDtime:1p3} - \ref{CRRtime:1p3} contain running time (in seconds) corresponding to Tables \ref{FDblad:1p3} - \ref{CRRblad:1p3}. Let us observe that the FD method is not significantly dependent on the scale parameter $\beta$ - that is, for each investigated value of $\beta$ the method returns the same values approximately at the same time. However, the MC and CRR methods depend significantly on $\beta$, where the sensitivity to this parameter increases for smaller $\alpha$. In particular, both methods work in a "long" time for the "small" values of $\beta$. Moreover, they are not precise for "big" $\beta$. 
		Optimization of $\beta$ with respect to the MC and CRR methods seems to be an interesting matter. We can conclude from this numerical example that the FD method has an advantage over the MC and CRR methods, which is the lack of sensitivity on the scale parameter $\beta$.
		\begin{table}[!ht]\fontsize{10}{12}\selectfont
			\begin{center}
				\begin{tabular}{|c| c c c c c|} 
					\hline
					\backslashbox{$\alpha$}{$\beta$} &   $10^{-8}$  &  $10^{-3}$&$1$&$10^{3}$&$10^{8}$\\	\hline
					$0.67$&   $3.81$&   $3.81$&   $3.81$ &  $3.81$&   $3.81$\\	
					$0.7$&   $3.68$&   $3.68$&   $3.68$ &  $3.68$&   $3.68$\\
					$0.8$&   $3.34$&   $3.34$&   $3.34$&   $3.34$&   $3.34$\\
					$0.9$&   $3.23$ &   $3.23$  &   $3.23$  &   $3.23$  &   $3.23$
					\\	
					$0.99$&   $3.39$ &   $3.39$  &   $3.39$  &   $3.39$  &   $3.39$
					\\		\hline
				\end{tabular}
				\caption{\label{FDblad:1p3}A relative error (multiplied by $10^{3}$) for the FD method in dependence on $\alpha$ and $\beta$.}
			\end{center}
		\end{table}	
		\begin{table}[!ht]\fontsize{10}{12}\selectfont
			\begin{center}
				\begin{tabular}{|c| c c c c c|} 
					\hline
					\backslashbox{$\alpha$}{$\beta$} &   $10^{-8}$  &  $10^{-3}$&$1$&$10^{3}$&$10^{8}$\\	\hline
					$0.67$& $8.67$&  $1.71$&  $3.87$&  $13.17$&  $440.75$\\
					$0.7$&   $2.34$&  $0.39$&  $4.96$&  $10.97$&  $436.89$\\
					$0.8$&  $2.85$&  $1.71$&  $9.29$&  $4.76$&  $89.9$ \\
					$0.9$&   $6.32$&  $2.61$&  $2.3$&  $8.08$&  $16.52$
					\\
					$0.99$&   $5.09$&  $5.62$ & $5.49$ & $4.81$&  $5.39$
					\\		\hline			
				\end{tabular}
				\caption{\label{MCblad:1p3}A relative error (multiplied by $10^{3}$) for the MC method in dependence on $\alpha$ and $\beta$.}
			\end{center}
		\end{table}
		\begin{table}[!ht]\fontsize{10}{12}\selectfont
			\begin{center}
				\begin{tabular}{|c| c c c c c|} 
					\hline
					\backslashbox{$\alpha$}{$\beta$} &   $10^{-8}$  &  $10^{-3}$&$1$&$10^{3}$&$10^{8}$\\	\hline	
					$0.67$&   $2.94$&  $0.35$&  $5.83$&  $13.96$&  $440.75$\\
					$0.7$&   $0.74$&  $4.41$&  $3.02$&  $18.49$&  $436.89$\\
					$0.8$&  $0.65$&  $2.49$&  $2.44$&  $7.53$&  $86.73$
					\\
					$0.9$&   $2.37$&  $1.74$&  $1.38$&  $6.05$&  $13.03$
					\\	
					$0.99$&   $2.74$&  $3.9$&  $3.45$&  $4.97$&  $4.32$
					\\		\hline
				\end{tabular}
				\caption{\label{CRRblad:1p3}A relative error (multiplied by $10^{3}$) for the CRR method in dependence on $\alpha$ and $\beta$.}
			\end{center}
		\end{table}
		\begin{table}[!ht]\fontsize{10}{12}\selectfont
			\begin{center}
				\begin{tabular}{|c| c c c c c|} 
					\hline
					\backslashbox{$\alpha$}{$\beta$} &   $10^{-8}$  &  $10^{-3}$&$1$&$10^{3}$&$10^{8}$\\	\hline	
					$0.67$&   $0.71$&   $0.63$&   $0.64$&   $0.66$&   $0.66$\\
					$0.7$&   $1.43$&   $1.44$&   $1.38$&   $1.49$&   $1.4$\\
					$0.8$&  $0.65$&   $0.56$&   $0.45$&   $0.46$&   $0.55$
					\\
					$0.9$&   $0.84$&   $0.77$&   $0.94$&   $0.91$&   $0.81$
					\\	
					$0.99$&   $0.6$ &  $0.63$&   $0.78$&   $   0.61$&   $0.63$
					\\		\hline
				\end{tabular}
				\caption{\label{FDtime:1p3}Computation time of FD metod corresponding to Table \ref{FDblad:1p3}.}
			\end{center}
		\end{table}
		\newpage\begin{table}[!ht]\fontsize{10}{12}\selectfont
			\begin{center}
				\begin{tabular}{|c| c c c c c|} 
					\hline
					\backslashbox{$\alpha$}{$\beta$} &   $10^{-8}$  &  $10^{-3}$&$1$&$10^{3}$&$10^{8}$\\	\hline	
					$0.67$&	$1971.66$&$   34.2$&$   3.27 $&$  0.44 $&$ 120.21$\\
					$0.7$&  $1127.43$&$   27.51$&$   3.4 $&$  0.55 $&$  2.64$\\
					$0.8$&  $127.52$&$   11.24$&$   2.8$&$   0.82$&$   0.23$
					\\
					$0.9$&   $14.37$&$   4.28$&$   2.27$&$   1.15$&$  0.5$
					\\	
					$0.99$&   $2.36$&$   2.1$&$   1.98$&$   1.98 $&$ 1.78$
					\\		\hline
				\end{tabular}
				\caption{\label{MCtime:1p3}Computation time of MC metod corresponding to Table \ref{MCblad:1p3}.}
			\end{center}
		\end{table}
		\begin{table}[!ht]\fontsize{10}{12}\selectfont
			\begin{center}
				\begin{tabular}{|c| c c c c c|} 
					\hline
					\backslashbox{$\alpha$}{$\beta$} &   $10^{-8}$  &  $10^{-3}$&$1$&$10^{3}$&$10^{8}$\\	\hline	
					$0.67$&   $2792.08$&$   36.98$&$   4.19 $&$  0.76 $&$  140.32$\\
					$0.7$&   $1813.15$&  $62.07$&  $ 8.51$&  $   1.77$&  $  6.65$\\
					$0.8$&  $128.61$&  $12.49$&  $3.98$&  $0.99$&  $0.44$
					\\
					$0.9$&   $12.89$&  $ 12.04$&  $ 6.23$&  $ 1.33$&  $ 0.64$
					\\	
					$0.99$&   $2.4$ &  $2.18$  & $2.06$  & $1.94$  & $1.75$
					\\		\hline
				\end{tabular}
				\caption{\label{CRRtime:1p3}Computation time of CRR metod corresponding to Table \ref{CRRblad:1p3}.}
			\end{center}
		\end{table}
	\end{Example}	

\newpage
\section{Summary}
In this paper:
\begin{itemize}
	\item[--] We have shown that the solution of the fractional tempered B-S equation is equal to the fair price of the European option with respect to $\mathbb{Q}$ in the tsB-S model.
	\item[--] We have introduced a weighted numerical scheme for this equation. It allows us to approximate the fair price of European call options in the tsB-S model.
	\item[--] We have given conditions under which the discrete scheme is stable and convergent.
	\item[--] We have discussed the problem of the optimal discretization parameter $\theta$ in dependence of the subdiffusion parameter $\alpha$ with tempering parameter  $\lambda$.
	\item[--] We have presented some numerical examples to illustrate the introduced theory.
\end{itemize}
We believe that the numerical techniques presented in this paper can be successfully repeated for other fractional diffusion-type problems.
\section*{Acknowledgments}
This research was partially supported by NCN Sonata Bis 9 grant nr 2019/34/E/ST1/00360.

\bibliographystyle{plain}
\bibliography{artykulTempSubDrugaP}

\begin{thebibliography}{10}

\bibitem{Meser}
M.~Alrawashdeh, J.~Kelly, M.~Meerschaert, and H.~Scheffler.
\newblock Applications of inverse tempered stable subordinators.
\newblock {\em Computers \& Mathematics with Applications}, 73(6):892--905,
  2017.

\bibitem{angelini2014delta}
F.~Angelini and S.~Herzel.
\newblock Delta hedging in discrete time under stochastic interest rate.
\newblock {\em Journal of Computational and Applied Mathematics}, 259:385--393,
  2014.

\bibitem{arregui2012numerical}
I.~Arregui and C.~V{\'a}zquez.
\newblock Numerical solution of an optimal investment problem with proportional
  transaction costs.
\newblock {\em Journal of Computational and Applied Mathematics},
  236(12):2923--2937, 2012.

\bibitem{jaCRR}
M.~Balcerek, G.~Krzy{\.z}anowski, and M.~Magdziarz.
\newblock About subordinated generalizations of 3 classical models of option
  pricing, 2021.

\bibitem{boen2020european}
L.~Boen.
\newblock European rainbow option values under the two-asset {Merton}
  jump-diffusion model.
\newblock {\em Journal of Computational and Applied Mathematics}, 364:112344,
  2020.

\bibitem{bor}
S.~Borak, A.~Misiorek, and R.~Weron.
\newblock Models for heavy-tailed asset returns.
\newblock In {\em Statistical tools for finance and insurance}, pages 21--55.
  Springer, 2011.

\bibitem{4}
R.~Cont and P.~Tankov.
\newblock {Financial Modelling with Jump Processes, Chapman \& Hall/CRC
  Financ}.
\newblock {\em Math. Ser}, 2004.

\bibitem{costabile2014option}
M.~Costabile, A.~Leccadito, I.~Massab{\'o}, and E.~Russo.
\newblock Option pricing under regime-switching jump--diffusion models.
\newblock {\em Journal of Computational and Applied Mathematics}, 256:152--167,
  2014.

\bibitem{de2021pseudospectral}
J.~de~Frutos and V.~Gat{\'o}n.
\newblock A pseudospectral method for option pricing with transaction costs
  under exponential utility.
\newblock {\em Journal of Computational and Applied Mathematics}, 394:113541,
  2021.

\bibitem{eliazar}
I.~Eliazar and J.~Klafter.
\newblock Spatial gliding, temporal trapping, and anomalous transport.
\newblock {\em Physica D: Nonlinear Phenomena}, 187(1-4):30--50, 2004.

\bibitem{elliott2014pricing}
R.~Elliott, T.~Siu, and L.~Chan.
\newblock On pricing barrier options with regime switching.
\newblock {\em Journal of Computational and Applied Mathematics}, 256:196--210,
  2014.

\bibitem{fama}
E.~Fama.
\newblock Risk, return and equilibrium: some clarifying comments.
\newblock {\em The Journal of Finance}, 23(1):29--40, 1968.

\bibitem{fia}
FIA.
\newblock {Global} {Futures} and {Options} {Trading} {Reaches} {Record} {Level}
  in 2019, 2019.

\bibitem{16}
V.~Gonchar, A.~Chechkin, E.~Sorokovoi, et~al.
\newblock Stable l{\'e}vy distributions of the density and potential
  fluctuations in the edge plasma of the {U-3M} torsatron.
\newblock {\em Plasma Physics Reports}, 29(5):380--390, 2003.

\bibitem{goncu2014efficient}
A.~G{\"o}nc{\"u} and G.~{\"O}kten.
\newblock Efficient simulation of a multi-factor stochastic volatility model.
\newblock {\em Journal of Computational and Applied Mathematics}, 259:329--335,
  2014.

\bibitem{heston2000rate}
S.~Heston and G.~Zhou.
\newblock On the rate of convergence of discrete-time contingent claims.
\newblock {\em Mathematical Finance}, 10(1):53--75, 2000.

\bibitem{janczura}
J.~Janczura and A.~Wy{\l}oma{\'n}ska.
\newblock Subdynamics of financial data from fractional {Fokker-Planck}
  equation.
\newblock {\em Acta Physica Polonica B}, 40(5):1341--1351, 2009.

\bibitem{15}
A.~Janicki and A.~Weron.
\newblock Can one see $\alpha$-stable variables and processes?
\newblock {\em Statistical Science}, pages 109--126, 1994.

\bibitem{ken1999levy}
S.~Ken-Iti.
\newblock {\em L{\'e}vy processes and infinitely divisible distributions}.
\newblock Cambridge University Press, 1999.

\bibitem{kilbas2006theory}
A.~Kilbas et~al.
\newblock {\em Theory and applications of fractional differential equations},
  volume 204.
\newblock Elsevier, 2006.

\bibitem{koleva2017numerical}
M.~Koleva and L.~Vulkov.
\newblock A numerical study for optimal portfolio regime-switching model i. 2d
  {Black--Scholes} equation with an exponential non-linear term.
\newblock {\em Journal of Computational and Applied Mathematics}, 318:538--549,
  2017.

\bibitem{ja}
G.~Krzy{\.z}anowski.
\newblock Selected applications of differential equations in {Vanilla Options}
  valuation.
\newblock {\em Mathematica Applicanda}, 46(2), 2018.

\bibitem{ja3}
G.~Krzy{\.z}anowski and M.~Magdziarz.
\newblock A computational weighted finite difference method for {American} and
  barrier options in subdiffusive {Black--Scholes} model.
\newblock {\em Communications in Nonlinear Science and Numerical Simulation},
  96:105676, 2020.

\bibitem{ja2}
G.~Krzy{\.z}anowski, M.~Magdziarz, and {\L}.~P{\l}ociniczak.
\newblock A weighted finite difference method for subdiffusive {Black--Scholes}
  model.
\newblock {\em Computers \& Mathematics with Applications}, 80(5):653--670,
  2020.

\bibitem{ladde2009development}
G.~Ladde and L.~Wu.
\newblock {Development} of modified geometric {Brownian} motion models by using
  stock price data and basic statistics.
\newblock {\em Nonlinear Analysis: Theory, Methods \& Applications},
  71(12):e1203--e1208, 2009.

\bibitem{Lin}
Y.~Lin and C.~Xu.
\newblock Finite difference/spectral approximations for the time-fractional
  diffusion equation.
\newblock {\em Journal of computational physics}, 225(2):1533--1552, 2007.

\bibitem{MM}
M.~Magdziarz.
\newblock {Black-Scholes} formula in subdiffusive regime.
\newblock {\em Journal of Statistical Physics}, 136(3):553--564, 2009.

\bibitem{Gajda}
M.~Magdziarz and J.~Gajda.
\newblock Anomalous dynamics of {Black-Scholes} model time-changed by inverse
  subordinators.
\newblock {\em Acta Physica Polonica B}, 43(5), 2012.

\bibitem{malhotra2022pricing}
G.~Malhotra, R.~Srivastava, and H.~Taneja.
\newblock Pricing of the geometric asian options under a multifactor stochastic
  volatility model.
\newblock {\em Journal of Computational and Applied Mathematics}, 406:113986,
  2022.

\bibitem{mandelbrot}
B.~Mandelbrot.
\newblock The variation of certain speculative prices.
\newblock In {\em Fractals and scaling in finance}, pages 371--418. Springer,
  1997.

\bibitem{McDonald}
R.~McDonald, M.~Cassano, and R.~Fahlenbrach.
\newblock {\em Derivatives markets}.
\newblock Addison-Wesley Boston, 2006.

\bibitem{17}
T.~Mizuuchi, V.~Chechkin, et~al.
\newblock Edge fluctuation studies in {Heliotron J}.
\newblock {\em Journal of nuclear materials}, 337:332--336, 2005.

\bibitem{orzel}
S.~Orze{\l} and A.~Weron.
\newblock Calibration of the subdiffusive {Black-Scholes} model.
\newblock {\em Acta Phys. Pol. B}, 41(5):1051--1059, 2010.

\bibitem{Wyl}
S.~Orze{\l} and A.~Wy{\l}oma{\'n}ska.
\newblock Calibration of the subdiffusive arithmetic {Brownian} motion with
  tempered stable waiting-times.
\newblock {\em Journal of Statistical Physics}, 143(3):447, 2011.

\bibitem{rachev1}
S.~Rachev and Frank~J. Menn, C.
\newblock {\em Fat-tailed and skewed asset return distributions: implications
  for risk management, portfolio selection, and option pricing}, volume 139.
\newblock John Wiley \& Sons, 2005.

\bibitem{rachev2}
S.~Rachev and S.~Mittnik.
\newblock {\em Stable {Paretian} models in finance}, volume~7.
\newblock Wiley, 2000.

\bibitem{rosinski}
J.~Rosi{\'n}ski.
\newblock Tempering stable processes.
\newblock {\em Stochastic processes and their applications}, 117(6):677--707,
  2007.

\bibitem{stanislavsky2003black}
A~Stanislavsky.
\newblock {Black--Scholes} model under subordination.
\newblock {\em Physica A: Statistical Mechanics and its Applications},
  318(3-4):469--474, 2003.

\bibitem{market}
A.~Stankovska.
\newblock {Global Derivatives Market}.
\newblock {\em SEEU Review}, 12, 01 2016.

\bibitem{18}
B.~Stuck and B.~Kleiner.
\newblock A statistical analysis of telephone noise.
\newblock {\em Bell System Technical Journal}, 53(7):1263--1320, 1974.

\bibitem{wang2021efficient}
W.~Wang, M.~Mao, and Z.~Wang.
\newblock An efficient variable step-size method for options pricing under
  jump-diffusion models with nonsmooth payoff function.
\newblock {\em ESAIM: Mathematical Modelling and Numerical Analysis},
  55(3):913--938, 2021.

\bibitem{wilmott1995mathematics}
P.~Wilmott, S.~Howison, and J.~Dewynne.
\newblock {\em The mathematics of financial derivatives: a student
  introduction}.
\newblock Cambridge university press, 1995.

\bibitem{yan2021utility}
D.~Yan and X.~Lu.
\newblock Utility-indifference pricing of european options with proportional
  transaction costs.
\newblock {\em Journal of Computational and Applied Mathematics}, 397:113639,
  2021.

\bibitem{zhang}
H.~Zhang, F.~Liu, I.~Turner, and Q.~Yang.
\newblock Numerical solution of the time fractional {Black--Scholes} model
  governing {European} options.
\newblock {\em Computers \& Mathematics with Applications}, 71(9):1772--1783,
  2016.

\end{thebibliography}

\end{document}